%% file: art2.tex
\documentstyle[12pt,amssymbo]{article}
\newenvironment{Proof}%
   {{\sc Proof: }}%
   { \hfill\cube\par}

   \font\Fraktur=eufm10 scaled\magstep1          % for display- and textstyle
   \newcommand{\fraktur}[1]{\mbox{\Fraktur #1}}  %
   \font\Fraktu=eufm7 scaled\magstep1            % for scriptstyle
   \newcommand{\fraktu}[1]{\mbox{\Fraktu #1}}    %
   \font\Frakt=eufm5 scaled\magstep1             % for scriptscriptstyle
   \newcommand{\frakt}[1]{\mbox{\Frakt #1}}      %
   \def\fr#1{\mathchoice{\fraktur {#1}}            % displaystyle
                        {\fraktur {#1}}            % textstyle
                        {\fraktu {#1}}             % scriptstyle
                        {\frakt {#1}}  }           % scriptscriptstyle

\newcommand{\Hom}{\mbox{\rm Hom}}
\newcommand{\Bi}{\mbox{\rm Im}}
\newcommand{\mapsfrom}{\hbox{$\leftarrow$\kern -0.15em \vrule width0.8pt height0.97ex depth-0.07ex}}

\newcommand{\pr}{\mbox{\rm pr}}
\newcommand{\plex}{^\bullet}
\newcommand{\zwind}[2]{\vbox{\baselineskip 0pt\hbox{$\scriptstyle #1$}\hbox{$\scriptstyle #2$}}}
\newtheorem{Definition}{Definition}[section]
\newtheorem{Theorem and Definition}[Definition]{Theorem and Definition}
\newtheorem{Definition and Lemma}[Definition]{Definition and Lemma}
\newtheorem{Example}[Definition]{Example}
\newtheorem{Theorem}[Definition]{Theorem}
\newtheorem{Proposition}[Definition]{Proposition}
\newtheorem{Remark}[Definition]{Remark}
\newtheorem{Lemma}[Definition]{Lemma}
\newtheorem{Corollary}[Definition]{Corollary}
\newcommand{\id}{\mbox{\it id}}
\newcommand{\const}{\mbox{\rm const}}
\newcommand{\si}[1]{(#1)}

\newcommand{\skp}[2]{\langle #1,#2 \rangle}
\newcommand{\cube}{{$\sqcap$\kern -0.7em$\sqcup$}}
\newcommand{\ou}[1]{\mathchoice%
  {\vtop{\baselineskip 0pt \hbox{\hfil$\,\otimes\,$\hfil}\hbox{\hfil$\scriptstyle #1$\hfil}%
         \vskip\lineskip}}%
  {\vtop{\baselineskip 0pt \hbox{\hfil$\,\otimes\,$\hfil}\hbox{\hfil$\scriptstyle #1$\hfil}%
         \vskip\lineskip}}%
  {\vtop{\baselineskip 0pt \hbox{$\scriptstyle\otimes$}%
         \hbox{$\scriptscriptstyle #1$}\vskip\lineskip}}%        
  {\vtop{\baselineskip 0pt \hbox{$\scriptscriptstyle\otimes$}%
         \hbox{$\scriptscriptstyle #1$}\vskip\lineskip}}}
\newcommand{\co}[1]{\mathchoice%
  {\vtop{\baselineskip 0pt \hbox{\hfil\,\cube\,\hfil}\hbox{\hfil$\scriptstyle #1$\hfil}%
         \vskip\lineskip}}%
  {\vtop{\baselineskip 0pt \hbox{\hfil\,\cube\,\hfil}\hbox{\hfil$\scriptstyle #1$\hfil}%
         \vskip\lineskip}}%
  {\vtop{\baselineskip 0pt \hbox{\scriptsize\cube}%
         \hbox{$\scriptscriptstyle #1$}\vskip\lineskip}}%        
  {\vtop{\baselineskip 0pt \hbox{\tiny\cube}%
         \hbox{$\scriptscriptstyle #1$}\vskip\lineskip}}}

\newcommand{\bpro}{{\noindent \sc Proof}: }
\newcommand{\epro}{\hfill $ \Box $ \parskip0.5cm \par \noindent \parskip1.0ex }
\newcommand{\bdm}{\begin{displaymath}}
\newcommand{\edm}{\end{displaymath} \noindent}
\newcommand{\be}{\begin{equation}}
\newcommand{\ee}{\end{equation} \noindent}
\newcommand{\bea}{\begin{eqnarray}}
\newcommand{\eea}{\end{eqnarray} \noindent}
\newcommand{\beas}{\begin{eqnarray*}}
\newcommand{\eeas}{\end{eqnarray*} \noindent}
\newcommand{\bmd}{\begin{eqnarray*} \begin{array}{lcclcl} }
\newcommand{\emd}{\end{array} \end{eqnarray*} \noindent}
\newcommand{\ot}{\otimes}
\newcommand{\bThe}{\begin{Theorem}}
\newcommand{\eThe}{\end{Theorem}}
\newcommand{\bPro}{\begin{Proposition}}
\newcommand{\ePro}{\end{Proposition}}
\newcommand{\bLem}{\begin{Lemma}}
\newcommand{\eLem}{\end{Lemma}}
\newcommand{\bCor}{\begin{Corollary}}
\newcommand{\eCor}{\end{Corollary}}
\newcommand{\bRem}{\begin{Remark}}
\newcommand{\eRem}{\end{Remark}}
\newcommand{\bExa}{\begin{Example}}
\newcommand{\eExa}{\end{Example}}
\newcommand{\bThD}{\begin{Theorem and Definition}}
\newcommand{\eThD}{\end{Theorem and Definition}}
\newcommand{\bDef}{\begin{Definition}}
\newcommand{\eDef}{\end{Definition}}
\newcommand{\Hopf}{$^H_H$Hopf$\,^H_H$-module}
\newcommand{\marg}[1]{}

\newcommand{\N}{{\Bbb{N}}}

\newcommand{\h}{\fr h}

\newcommand{\rcofix}[2]{{#1}^{\mbox{\rm\scriptsize co }{#2}}}
\newcommand{\lcofix}[2]{{^{\mbox{\rm\scriptsize co }{#2}}}{#1}}
\newcommand{\HMod}[4]{{^{#1}_{#3}}{\cal M}^{#2}_{#4}}

\input{catmac}
\newcommand{\nabb}[1]{%
   \horsize{\tempcounta}{{\scriptstyle{#1}}}%
   \xext = \tempcounta%
   \yext = 0%
   \horsize{\tempcountb}{{00}}%
   \advance \xext by\tempcountb%
   \adjust[\scriptstyle{#1}`\scriptstyle{#1};`;`;\scriptstyle{#1}`]%
   \begin{picture}(\xext,\yext)(\xoff,\yoff)%
     \puthmorphism(0,0)[``{\scriptstyle{#1}}]{\xext}1a%
   \end{picture}%
   }
\renewcommand{\epsilon}{\varepsilon}
\begin{document}
\include{title}
\tableofcontents
\include{intro}

\include{sect1}

\include{sect2}
\include{sect3}
\include{sect4}

\include{sect5}
\bibliographystyle{amsplain}
\bibliography{Algebra,mathlib}
\end{document}

%% file: catmac.tex
\long\def\ig#1{\relax}
\ig{Thanks to Roberto Minio for this def'n.  Compare the def'n of
\comment in AMSTeX.}
 
\newcount \coefa
\newcount \coefb
\newcount \coefc
\newdimen\tempdimen
\newdimen\xlen
\newdimen\ylen
 
\newcount\tempcounta
\newcount\tempcountb
\newcount\tempcountc
\newcount\tempcountd
\newcount\tempcounte
\newcount\tempcountf
\newcount\xext
\newcount\yext
\newcount\xoff
\newcount\yoff
\newcount\gap%
\newcount\arrowtypea
\newcount\arrowtypeb
\newcount\arrowtypec
\newcount\arrowtyped
\newcount\arrowtypee
\newcount\height
\newcount\width
\newcount\xpos
\newcount\ypos
\newcount\run
\newcount\rise
\newcount\arrowlength
\newcount\halflength
\newcount\arrowtype
\newsavebox{\tempboxa}%
\newsavebox{\tempboxb}%
\newsavebox{\tempboxc}%

\catcode`@=11 % unhide @sign definitions%
\def\settoheight#1#2{\setbox\@tempboxa\hbox{#2}#1\ht\@tempboxa\relax}%
\def\settodepth#1#2{\setbox\@tempboxa\hbox{#2}#1\dp\@tempboxa\relax}%
\let\ifnextchar=\@ifnextchar
\catcode`@=12 %@signs are no longer letters%
 
\def\putbox(#1,#2)#3{\put(#1,#2){\makebox(0,0){#3}}}

\def\setsqparms[#1`#2`#3`#4;#5`#6]{%
\settripairparms[#1`#2`#3`#4`1;#6]%
\width #5
}
 
\def\settriparms[#1`#2`#3;#4]{\settripairparms[#1`#2`#3`1`1;#4]}%

\def\settripairparms[#1`#2`#3`#4`#5;#6]{%
\arrowtypea #1
\arrowtypeb #2
\arrowtypec #3
\arrowtyped #4
\arrowtypee #5
\height #6
\width #6
}
 
\def\mvector(#1,#2)#3{%%
\put(0,0){\vector(#1,#2){#3}}%
\put(0,0){\vector(#1,#2){30}}%
}
\def\evector(#1,#2)#3{{%%
\arrowlength #3
\put(0,0){\vector(#1,#2){\arrowlength}}%
\advance \arrowlength by-30
\put(0,0){\vector(#1,#2){\arrowlength}}%
}}

\def\horsize#1#2{%
\settowidth{\tempdimen}{$#2$}%
#1=\tempdimen
\divide #1 by\unitlength
}
 
\def\vertsize#1#2{%
\settoheight{\tempdimen}{$#2$}%
#1=\tempdimen
\settodepth{\tempdimen}{$#2$}%
\advance #1 by\tempdimen
\divide #1 by\unitlength
}

\def\vertadjust[#1`#2`#3]{%
\vertsize{\tempcounta}{#1}%
\vertsize{\tempcountb}{#2}%
\ifnum \tempcounta<\tempcountb \tempcounta=\tempcountb \fi
\divide\tempcounta by2
\vertsize{\tempcountb}{#3}%
\ifnum \tempcountb>0 \advance \tempcountb by20 \fi
\ifnum \tempcounta<\tempcountb \tempcounta=\tempcountb \fi
}
 
\def\horadjust[#1`#2`#3]{%
\horsize{\tempcounta}{#1}%
\horsize{\tempcountb}{#2}%
\ifnum \tempcounta<\tempcountb \tempcounta=\tempcountb \fi
\divide\tempcounta by20
\horsize{\tempcountb}{#3}%
\ifnum \tempcountb>0 \advance \tempcountb by60 \fi
\ifnum \tempcounta<\tempcountb \tempcounta=\tempcountb \fi
}
 
\ig{ In this procedure, #1 is the paramater that sticks out all the way,
#2 sticks out the least and #3 is a label sticking out half way.  #4 is
the amount of the offset.}
 
\def\sladjust[#1`#2`#3]#4{%
\tempcountc=#4
\horsize{\tempcounta}{#1}%
\divide \tempcounta by2
\horsize{\tempcountb}{#2}%
\divide \tempcountb by2
\advance \tempcountb by-\tempcountc
\ifnum \tempcounta<\tempcountb \tempcounta=\tempcountb\fi
\divide \tempcountc by2
\horsize{\tempcountb}{#3}%
\advance \tempcountb by-\tempcountc
\ifnum \tempcountb>0 \advance \tempcountb by80\fi
\ifnum \tempcounta<\tempcountb \tempcounta=\tempcountb\fi
\advance\tempcounta by20
}
 
\def\putvector(#1,#2)(#3,#4)#5#6{{%
\xpos=#1
\ypos=#2
\run=#3
\rise=#4
\arrowlength=#5
\arrowtype=#6
\ifnum \arrowtype<0
    \ifnum \run=0
        \advance \ypos by-\arrowlength
    \else
        \tempcounta \arrowlength
        \multiply \tempcounta by\rise
        \divide \tempcounta by\run
        \ifnum\run>0
            \advance \xpos by\arrowlength
            \advance \ypos by\tempcounta
        \else
            \advance \xpos by-\arrowlength
            \advance \ypos by-\tempcounta
        \fi
    \fi
    \multiply \arrowtype by-1
    \multiply \rise by-1
    \multiply \run by-1
\fi
\ifnum \arrowtype=1
    \put(\xpos,\ypos){\vector(\run,\rise){\arrowlength}}%
\else\ifnum \arrowtype=2
    \put(\xpos,\ypos){\mvector(\run,\rise)\arrowlength}%
\else\ifnum\arrowtype=3
    \put(\xpos,\ypos){\evector(\run,\rise){\arrowlength}}%
\fi\fi\fi
}}
 
\def\bfig{\begin{picture}(\xext,\yext)(\xoff,\yoff)}
\def\efig{\end{picture}}

\def\putsplitvector(#1,#2)#3#4{%%
\xpos #1
\ypos #2
\arrowtype #4
\halflength #3
\arrowlength #3
\gap 140
\advance \halflength by-\gap
\divide \halflength by2
\ifnum \arrowtype=1
    \put(\xpos,\ypos){\line(0,-1){\halflength}}%
    \advance\ypos by-\halflength
    \advance\ypos by-\gap
    \put(\xpos,\ypos){\vector(0,-1){\halflength}}%
\else\ifnum \arrowtype=2
    \put(\xpos,\ypos){\line(0,-1)\halflength}%
    \put(\xpos,\ypos){\vector(0,-1)3}%
    \advance\ypos by-\halflength
    \advance\ypos by-\gap
    \put(\xpos,\ypos){\vector(0,-1){\halflength}}%
\else\ifnum\arrowtype=3
    \put(\xpos,\ypos){\line(0,-1)\halflength}%
    \advance\ypos by-\halflength
    \advance\ypos by-\gap
    \put(\xpos,\ypos){\evector(0,-1){\halflength}}%
\else\ifnum \arrowtype=-1
    \advance \ypos by-\arrowlength
    \put(\xpos,\ypos){\line(0,1){\halflength}}%
    \advance\ypos by\halflength
    \advance\ypos by\gap
    \put(\xpos,\ypos){\vector(0,1){\halflength}}%
\else\ifnum \arrowtype=-2
    \advance \ypos by-\arrowlength
    \put(\xpos,\ypos){\line(0,1)\halflength}%
    \put(\xpos,\ypos){\vector(0,1)3}%
    \advance\ypos by\halflength
    \advance\ypos by\gap
    \put(\xpos,\ypos){\vector(0,1){\halflength}}%
\else\ifnum\arrowtype=-3
    \advance \ypos by-\arrowlength
    \put(\xpos,\ypos){\line(0,1)\halflength}%
    \advance\ypos by\halflength
    \advance\ypos by\gap
    \put(\xpos,\ypos){\evector(0,1){\halflength}}%
\fi\fi\fi\fi\fi\fi
}
 
\def\setpos(#1,#2){\xpos=#1 \ypos#2}
 
\def\putmorphism(#1)(#2,#3)[#4`#5`#6]#7#8#9{{%
\run #2
\rise #3
\ifnum\rise=0
  \puthmorphism(#1)[#4`#5`#6]{#7}{#8}{#9}%
\else\ifnum\run=0
  \putvmorphism(#1)[#4`#5`#6]{#7}{#8}{#9}%
\else
\setpos(#1)%
\arrowlength #7
\arrowtype #8
\ifnum\run=0
\else\ifnum\rise=0
\else
\ifnum\run>0
    \coefa=1
\else
   \coefa=-1
\fi
\ifnum\arrowtype>0
   \coefb=0
   \coefc=-1
\else
   \coefb=\coefa
   \coefc=1
   \arrowtype=-\arrowtype
\fi
\width=2
\multiply \width by\run
\divide \width by\rise
\ifnum \width<0  \width=-\width\fi
\advance\width by60
\if l#9 \width=-\width\fi
\putbox(\xpos,\ypos){$#4$}%            %node 1
{\multiply \coefa by\arrowlength%      %node 2
\advance\xpos by\coefa
\multiply \coefa by\rise
\divide \coefa by\run
\advance \ypos by\coefa
\putbox(\xpos,\ypos){$#5$} }%
{\multiply \coefa by\arrowlength%      %label
\divide \coefa by2
\advance \xpos by\coefa
\advance \xpos by\width
\multiply \coefa by\rise
\divide \coefa by\run
\advance \ypos by\coefa
\if l#9%
   \put(\xpos,\ypos){\makebox(0,0)[r]{$#6$}}%
\else\if r#9%
   \put(\xpos,\ypos){\makebox(0,0)[l]{$#6$}}%
\fi\fi }%
{\multiply \rise by-\coefc%             %arrow
\multiply \run by-\coefc
\multiply \coefb by\arrowlength
\advance \xpos by\coefb
\multiply \coefb by\rise
\divide \coefb by\run
\advance \ypos by\coefb
\multiply \coefc by70
\advance \ypos by\coefc
\multiply \coefc by\run
\divide \coefc by\rise
\advance \xpos by\coefc
\multiply \coefa by140
\multiply \coefa by\run
\divide \coefa by\rise
\advance \arrowlength by\coefa
\ifnum \arrowtype=1
   \put(\xpos,\ypos){\vector(\run,\rise){\arrowlength}}%
\else\ifnum\arrowtype=2
   \put(\xpos,\ypos){\mvector(\run,\rise){\arrowlength}}%
\else\ifnum\arrowtype=3
   \put(\xpos,\ypos){\evector(\run,\rise){\arrowlength}}%
\fi\fi\fi}%
\fi\fi
\fi\fi}}

\def\puthmorphism(#1,#2)[#3`#4`#5]#6#7#8{{%
\xpos #1
\ypos #2
\width #6
\arrowlength #6
\putbox(\xpos,\ypos){$#3$\vphantom{$#4$}}%
{\advance \xpos by\arrowlength
\putbox(\xpos,\ypos){\vphantom{$#3$}$#4$}}%
\horsize{\tempcounta}{#3}%
\horsize{\tempcountb}{#4}%
\divide \tempcounta by2
\divide \tempcountb by2
\advance \tempcounta by30
\advance \tempcountb by30
\advance \xpos by\tempcounta
\advance \arrowlength by-\tempcounta
\advance \arrowlength by-\tempcountb
\putvector(\xpos,\ypos)(1,0){\arrowlength}{#7}%
\divide \arrowlength by2
\advance \xpos by\arrowlength
\vertsize{\tempcounta}{#5}%
\divide\tempcounta by2
\advance \tempcounta by20
\if a#8 %
   \advance \ypos by\tempcounta
   \put(\xpos,\ypos){\makebox(0,0){$#5$}}%
\else
   \advance \ypos by-\tempcounta
   \put(\xpos,\ypos){\makebox(0,0){$#5$}}%
\fi
}}
 
\def\putvmorphism(#1,#2)[#3`#4`#5]#6#7#8{{%
\xpos #1
\ypos #2
\arrowlength #6
\arrowtype #7
\settowidth{\xlen}{$#5$}%
\putbox(\xpos,\ypos){$#3$}%
{\advance \ypos by-\arrowlength
\putbox(\xpos,\ypos){$#4$}}%
{\advance\arrowlength by-140
\advance \ypos by-70
\ifdim\xlen>0pt
   \if m#8%
      \putsplitvector(\xpos,\ypos){\arrowlength}{\arrowtype}%
   \else
      \putvector(\xpos,\ypos)(0,-1){\arrowlength}{\arrowtype}%
   \fi
\else
   \putvector(\xpos,\ypos)(0,-1){\arrowlength}{\arrowtype}%
\fi}%
\ifdim\xlen>0pt
   \divide \arrowlength by2
   \advance\ypos by-\arrowlength
   \if l#8%
      \advance \xpos by-40
      \put(\xpos,\ypos){\makebox(0,0)[r]{$#5$}}%
   \else\if r#8%
      \advance \xpos by40
      \put(\xpos,\ypos){\makebox(0,0)[l]{$#5$}}%
   \else
      \putbox(\xpos,\ypos){$#5$}%
   \fi\fi
\fi
}}
 
\def\topadjust[#1`#2`#3]{%
\yoff=10
\vertadjust[#1`#2`{#3}]%
\advance \yext by\tempcounta
\advance \yext by 10
}
 
\def\botadjust[#1`#2`#3]{%
\vertadjust[#1`#2`{#3}]%
\advance \yext by\tempcounta
\advance \yoff by-\tempcounta
}
 
\def\leftadjust[#1`#2`#3]{%
\xoff=0
\horadjust[#1`#2`{#3}]%
\advance \xext by\tempcounta
\advance \xoff by-\tempcounta
}
 
\def\rightadjust[#1`#2`#3]{%
\horadjust[#1`#2`{#3}]%
\advance \xext by\tempcounta
}
 
\def\rightsladjust[#1`#2`#3]{%
\sladjust[#1`#2`{#3}]{\width}%
\advance \xext by\tempcounta
}
 
\def\leftsladjust[#1`#2`#3]{%
\xoff=0
\sladjust[#1`#2`{#3}]{\width}%
\advance \xext by\tempcounta
\advance \xoff by-\tempcounta
}
 
\def\adjust[#1`#2;#3`#4;#5`#6;#7`#8]{%
\topadjust[#1``{#2}]
\leftadjust[#3``{#4}]
\rightadjust[#5``{#6}]
\botadjust[#7``{#8}]}

\def\putsquare(#1)[#2`#3`#4`#5;#6`#7`#8`#9]{%
\setpos(#1)
\puthmorphism(\xpos,\ypos)[#4`#5`{#9}]{\width}{\arrowtyped}b%
\advance\ypos by \height
\puthmorphism(\xpos,\ypos)[#2`#3`{#6}]{\width}{\arrowtypea}a%
\putvmorphism(\xpos,\ypos)[``{#7}]{\height}{\arrowtypeb}l%
\advance\xpos by \width
\putvmorphism(\xpos,\ypos)[``{#8}]{\height}{\arrowtypec}r%
}
 
\def\square[#1`#2`#3`#4;#5`#6`#7`#8]{{%%            #5
\xext=\width                              %     #1------>#2
\yext=\height                             %      |       |
\topadjust[#1`#2`{#5}]%                   %      |       |
\botadjust[#3`#4`{#8}]%                   %    #6|       |#7
\leftadjust[#1`#3`{#6}]%                  %      |       |
\rightadjust[#2`#4`{#7}]%                 %      |       |
\begin{picture}(\xext,\yext)(\xoff,\yoff)%       v       v
\putsquare(0,0)[#1`#2`#3`#4;#5`#6`#7`{#8}]%     #3------>#4
\end{picture}%                                      #8
}}

\def\putptriangle(#1,#2)[#3`#4`#5;#6`#7`#8]{%
\xpos=#1 \ypos=#2
\advance\ypos by \height
\puthmorphism(\xpos,\ypos)[#3`#4`{#6}]{\height}{\arrowtypea}a%
\putvmorphism(\xpos,\ypos)[`#5`{#7}]{\height}{\arrowtypeb}l%
\advance\xpos by\height
\putmorphism(\xpos,\ypos)(-1,-1)[``{#8}]{\height}{\arrowtypec}r%
}
 
\def\ptriangle[#1`#2`#3;#4`#5`#6]{{%%            #4
\width=\height                         %      #1----->#2
\xext=\width                           %      |      /
\yext=\width                           %      |     /
\topadjust[#1`#2`{#4}]%                %    #5|    /#6
\botadjust[#3``]%                      %      |   /
\leftadjust[#1`#3`{#5}]%               %      |  /
\rightsladjust[#2`#3`{#6}]%            %      v v
\begin{picture}(\xext,\yext)(\xoff,\yoff)%    #3
\putptriangle(0,0)[#1`#2`#3;#4`#5`{#6}]%
\end{picture}%
}}

\def\putqtriangle(#1,#2)[#3`#4`#5;#6`#7`#8]{%
\xpos=#1 \ypos=#2
\advance\ypos by\height
\puthmorphism(\xpos,\ypos)[#3`#4`{#6}]{\height}{\arrowtypea}a%
\putmorphism(\xpos,\ypos)(1,-1)[``{#7}]{\height}{\arrowtypeb}l%
\advance\xpos by\height
\putvmorphism(\xpos,\ypos)[`#5`{#8}]{\height}{\arrowtypec}r%
}
 
\def\qtriangle[#1`#2`#3;#4`#5`#6]{{%%              #4
\width=\height                         %        #1----->#2
\xext=\width                           %         \      |
\yext=\height                          %          \     |
\topadjust[#1`#2`{#4}]%                %         #5\    |#6
\botadjust[#3``]%                      %            \   |
\leftsladjust[#1`#3`{#5}]%             %             \  |
\rightadjust[#2`#3`{#6}]%              %              v v
\begin{picture}(\xext,\yext)(\xoff,\yoff)%             #3
\putqtriangle(0,0)[#1`#2`#3;#4`#5`{#6}]%
\end{picture}%
}}

\def\putdtriangle(#1,#2)[#3`#4`#5;#6`#7`#8]{%
\xpos=#1 \ypos=#2
\puthmorphism(\xpos,\ypos)[#4`#5`{#8}]{\height}{\arrowtypec}b%
\advance\xpos by \height \advance\ypos by\height
\putmorphism(\xpos,\ypos)(-1,-1)[``{#6}]{\height}{\arrowtypea}l%
\putvmorphism(\xpos,\ypos)[#3``{#7}]{\height}{\arrowtypeb}r%
}
 
\def\dtriangle[#1`#2`#3;#4`#5`#6]{{%%                      #1
\width=\height                         %                  / |
\xext=\width                           %                 /  |
\yext=\height                          %              #4/   |#5
\topadjust[#1``]%                      %               /    |
\botadjust[#2`#3`{#6}]%                %              /     |
\leftsladjust[#2`#1`{#4}]%             %             v      v
\rightadjust[#1`#3`{#5}]%              %            #2----->#3
\begin{picture}(\xext,\yext)(\xoff,\yoff)%              #6
\putdtriangle(0,0)[#1`#2`#3;#4`#5`{#6}]%
\end{picture}%
}}

\def\putbtriangle(#1,#2)[#3`#4`#5;#6`#7`#8]{%
\xpos=#1 \ypos=#2
\puthmorphism(\xpos,\ypos)[#4`#5`{#8}]{\height}{\arrowtypec}b%
\advance\ypos by\height
\putmorphism(\xpos,\ypos)(1,-1)[``{#7}]{\height}{\arrowtypeb}r%
\putvmorphism(\xpos,\ypos)[#3``{#6}]{\height}{\arrowtypea}l%
}
 
\def\btriangle[#1`#2`#3;#4`#5`#6]{{%%                 #1
\width=\height                         %              | \
\xext=\width                           %              |  \
\yext=\height                          %            #4|   \#5
\topadjust[#1``]%                      %              |    \
\botadjust[#2`#3`{#6}]%                %              |     \
\leftadjust[#1`#2`{#4}]%               %              v      v
\rightsladjust[#3`#1`{#5}]%            %              #2----->#3
\begin{picture}(\xext,\yext)(\xoff,\yoff)%                #6
\putbtriangle(0,0)[#1`#2`#3;#4`#5`{#6}]%
\end{picture}%
}}

\def\putAtriangle(#1,#2)[#3`#4`#5;#6`#7`#8]{%
\xpos=#1 \ypos=#2
{\multiply \height by2
\puthmorphism(\xpos,\ypos)[#4`#5`{#8}]{\height}{\arrowtypec}b}%
\advance\xpos by\height \advance\ypos by\height
\putmorphism(\xpos,\ypos)(-1,-1)[#3``{#6}]{\height}{\arrowtypea}l%
\putmorphism(\xpos,\ypos)(1,-1)[``{#7}]{\height}{\arrowtypeb}r%
}
 
\def\Atriangle[#1`#2`#3;#4`#5`#6]{{%%             #1
\width=\height                         %         /   \
\xext=\width                           %        /     \
\yext=\height                          %     #4/       \#5
\topadjust[#1``]%                      %      /         \
\botadjust[#2`#3`{#6}]%                %     /           \
\multiply \xext by2 %                  %    v             v
\leftsladjust[#2`#1`{#4}]%             %   #2------------>#3
\rightsladjust[#3`#1`{#5}]%            %          #6
\begin{picture}(\xext,\yext)(\xoff,\yoff)%
\putAtriangle(0,0)[#1`#2`#3;#4`#5`{#6}]%
\end{picture}%
}}

\def\putAtrianglepair(#1,#2)[#3]{\xpos=#1 \ypos=#2%
\putAtrianglepairx[#3]}
\def\putAtrianglepairx[#1`#2`#3`#4;#5`#6`#7`#8`#9]{%
\puthmorphism(\xpos,\ypos)[#2`#3`{#8}]{\height}{\arrowtyped}b%
\advance\xpos by\height
\puthmorphism(\xpos,\ypos)[\phantom{#3}`#4`{#9}]{\height}{\arrowtypee}b%
\advance\ypos by\height
\putmorphism(\xpos,\ypos)(-1,-1)[#1``{#5}]{\height}{\arrowtypea}l%
\putvmorphism(\xpos,\ypos)[``{#6}]{\height}{\arrowtypeb}m%
\putmorphism(\xpos,\ypos)(1,-1)[``{#7}]{\height}{\arrowtypec}r%
}
 
\def\Atrianglepair[#1`#2`#3`#4;#5`#6`#7`#8`#9]{{%
\width=\height
\xext=\width
\yext=\height
\topadjust[#1``]%
\vertadjust[#2`#3`{#8}]%                      %            #1
\tempcountd=\tempcounta                       %           / | \
\vertadjust[#3`#4`{#9}]%                      %          /  |  \
\ifnum\tempcounta<\tempcountd                 %       #5/  #6   \#7
\tempcounta=\tempcountd\fi                    %        /    |    \
\advance \yext by\tempcounta                  %       /     |     \
\advance \yoff by-\tempcounta                 %      v      v      v
\multiply \xext by2 %                         %     #2----->#3----->#4
\leftsladjust[#2`#1`{#5}]%                    %         #8      #9
\rightsladjust[#4`#1`{#7}]%
\begin{picture}(\xext,\yext)(\xoff,\yoff)%
\putAtrianglepair(0,0)[#1`#2`#3`#4;#5`#6`#7`#8`{#9}]%
\end{picture}%
}}

\def\putVtriangle(#1,#2)[#3`#4`#5;#6`#7`#8]{%
\xpos=#1 \ypos=#2
\advance\ypos by\height
{\multiply\height by2
\puthmorphism(\xpos,\ypos)[#3`#4`{#6}]{\height}{\arrowtypea}a}%
\putmorphism(\xpos,\ypos)(1,-1)[`#5`{#7}]{\height}{\arrowtypeb}l%
\advance\xpos by\height
\advance\xpos by\height
\putmorphism(\xpos,\ypos)(-1,-1)[``{#8}]{\height}{\arrowtypec}r%
}
 
\def\Vtriangle[#1`#2`#3;#4`#5`#6]{{%%                  #4
\width=\height                         %        #1------------->#2
\xext=\width                           %         \             /
\yext=\height                          %          \           /
\topadjust[#1`#2`{#4}]%                %         #5\         /#6
\botadjust[#3``]%                      %            \       /
\multiply \xext by2 %                  %             \     /
\leftsladjust[#1`#3`{#5}]%             %              v   v
\rightsladjust[#2`#3`{#6}]%            %               #3
\begin{picture}(\xext,\yext)(\xoff,\yoff)%
\putVtriangle(0,0)[#1`#2`#3;#4`#5`{#6}]%
\end{picture}%
}}

\def\putVtrianglepair(#1,#2)[#3]{\xpos=#1 \ypos=#2%
\putVtrianglepairx[#3]}
\def\putVtrianglepairx[#1`#2`#3`#4;#5`#6`#7`#8`#9]{%
\advance\ypos by\height
\putmorphism(\xpos,\ypos)(1,-1)[`#4`{#7}]{\height}{\arrowtypec}l%
\puthmorphism(\xpos,\ypos)[#1`#2`{#5}]{\height}{\arrowtypea}a%
\advance\xpos by\height
\puthmorphism(\xpos,\ypos)[\phantom{#2}`#3`{#6}]{\height}{\arrowtypeb}a%
\putvmorphism(\xpos,\ypos)[``{#8}]{\height}{\arrowtyped}m%
\advance\xpos by\height
\putmorphism(\xpos,\ypos)(-1,-1)[``{#9}]{\height}{\arrowtypee}r%
}
 
\def\Vtrianglepair[#1`#2`#3`#4;#5`#6`#7`#8`#9]{{%
\xoff=0
\yoff=2 %                      %           #5      #6
\xext=\height                  %        #1----->#2----->#3
\width=\height                 %         \      |      /
\yext=\height                  %          \     |     /
\vertadjust[#1`#2`{#5}]%       %         #7\   #8    /#9
\tempcountd=\tempcounta        %            \   |   /
\vertadjust[#2`#3`{#6}]%       %             \  |  /
\ifnum\tempcounta<\tempcountd%                v v v
\tempcounta=\tempcountd\fi%                    #4
\advance \yext by\tempcounta
\botadjust[#4``]%
\multiply \xext by2
\leftsladjust[#1`#4`{#7}]%
\rightsladjust[#3`#4`{#9}]%
\begin{picture}(\xext,\yext)(\xoff,\yoff)%
\putVtrianglepair(0,0)[#1`#2`#3`#4;#5`#6`#7`#8`{#9}]%
\end{picture}%
}}

\def\putCtriangle(#1,#2)[#3`#4`#5;#6`#7`#8]{%
\xpos=#1 \ypos=#2
\advance\ypos by\height
\putmorphism(\xpos,\ypos)(1,-1)[``{#8}]{\height}{\arrowtypec}l%
\advance\xpos by\height
\advance\ypos by\height
\putmorphism(\xpos,\ypos)(-1,-1)[#3`#4`{#6}]{\height}{\arrowtypea}l%
{\multiply\height by 2
\putvmorphism(\xpos,\ypos)[`#5`{#7}]{\height}{\arrowtypeb}r}%
}
 
\def\Ctriangle[#1`#2`#3;#4`#5`#6]{{% %                    #1
\width=\height                          %                / |
\xext=\width                            %               /  |
\yext=\height                           %            #4/   |
\multiply \yext by2 %                   %             /    |
\topadjust[#1``]%                       %            /     |
\botadjust[#3``]%                       %           v      |
\sladjust[#2`#1`{#4}]{\width}%          %          #2      |#5
\tempcountd=\tempcounta                 %           \      |
\sladjust[#2`#3`{#6}]{\width}%          %            \     |
\ifnum \tempcounta<\tempcountd          %           #6\    |
\tempcounta=\tempcountd\fi              %              \   |
\advance \xext by\tempcounta            %               \  |
\advance \xoff by-\tempcounta           %                v v
\rightadjust[#1`#3`{#5}]%               %                 #3
\begin{picture}(\xext,\yext)(\xoff,\yoff)%
\putCtriangle(0,0)[#1`#2`#3;#4`#5`{#6}]%
\end{picture}%
}}

\def\putDtriangle(#1,#2)[#3`#4`#5;#6`#7`#8]{%
\xpos=#1 \ypos=#2
\advance\xpos by\height \advance\ypos by\height
\putmorphism(\xpos,\ypos)(-1,-1)[``{#8}]{\height}{\arrowtypec}r%
\advance\xpos by-\height \advance\ypos by\height
\putmorphism(\xpos,\ypos)(1,-1)[`#4`{#7}]{\height}{\arrowtypeb}r%
{\multiply\height by 2
\putvmorphism(\xpos,\ypos)[#3`#5`{#6}]{\height}{\arrowtypea}l}%
}
 
\def\Dtriangle[#1`#2`#3;#4`#5`#6]{{%%             #1
\width=\height                         %          | \
\xext=\height                          %          |  \
\yext=\height                          %          |   \#5
\multiply \yext by2 %                  %          |    \
\topadjust[#1``]%                      %          |     \
\botadjust[#3``]%                      %          |      v
\leftadjust[#1`#3`{#4}]%               %        #4|      #2
\sladjust[#2`#1`{#4}]{\height}%        %          |      /
\tempcountd=\tempcountd                %          |     /
\sladjust[#2`#3`{#6}]{\height}%        %          |    /#6
\ifnum \tempcounta<\tempcountd         %          |   /
\tempcounta=\tempcountd\fi             %          |  /
\advance \xext by\tempcounta           %          v v
\begin{picture}(\xext,\yext)(\xoff,\yoff)%        #3
\putDtriangle(0,0)[#1`#2`#3;#4`#5`{#6}]%
\end{picture}%
}}

\def\setrecparms[#1`#2]{\width=#1 \height=#2}%
\def\recurse[#1`#2`#3`#4;#5`#6`#7`#8`#9]{{%
\settowidth{\tempdimen}{#1}
\ifdim\tempdimen=0pt
  \savebox{\tempboxa}{\hbox{#2}}%
  \savebox{\tempboxb}{\hbox{#4}}%
  \savebox{\tempboxc}{\hbox{#7}}%
\else
  \savebox{\tempboxa}{\hbox{$\hbox{#1}\times\hbox{#2}$}}%
  \savebox{\tempboxb}{\hbox{$\hbox{#1}\times\hbox{#4}$}}%
  \savebox{\tempboxc}{\hbox{$\hbox{#1}\times\hbox{#7}$}}%
\fi
\tempcounte=\height
\divide\tempcounte by 2
\tempcountf=\tempcounte
\advance\tempcountf by \width
\xext=\tempcountf \yext=\height
\topadjust[#2`\usebox{\tempboxa}`{#5}]%
\botadjust[#4`\usebox{\tempboxb}`{#9}]%
\sladjust[#3`#2`{#6}]{\tempcounte}%
\tempcountd=\tempcounta
\sladjust[#3`#4`{#8}]{\tempcounte}%
\ifnum \tempcounta<\tempcountd
\tempcounta=\tempcountd\fi
\advance \xext by\tempcounta
\advance \xoff by-\tempcounta
\rightadjust[\usebox{\tempboxa}`\usebox{\tempboxb}`\usebox{\tempboxc}]%
\bfig
{\settriparms[-1`1`1;\tempcounte]%
\putCtriangle(0,0)[`#3`;#6`#7`{#8}]}%
\arrowtypea=-1 \arrowtypeb=0 \arrowtypec=1 \arrowtyped=-1
\putsquare(\tempcounte,0)[#2`\usebox{\tempboxa}`#4`\usebox{\tempboxb};%
#5``\usebox{\tempboxc}`#9]%
\efig
}}

%% file: title.tex
\begin{titlepage}
\title{Differential calculi on noncommutative bundles}
\author {Markus J.~Pflaum, Peter Schauenburg\\
  pflaum@rz.mathematik.uni-muenchen.de,\\
  schauen@rz.mathematik.uni-muenchen.de}
\date{Mathematisches Institut der Universit\"{a}t M\"{u}nchen\\
Theresienstra{\ss}e 39\\
80333 M\"{u}nchen 2, Germany\\
29.~Juli 1994}
\maketitle
\thispagestyle{empty}
\begin{abstract}
We introduce a category of noncommutative bundles.
To establish  geometry in this category
we construct suitable noncommutative differential calculi on these
bundles and study their
basic properties. 
Furthermore we define the notion of a connection with respect to a 
differential calculus
and consider questions of existence and uniqueness. 
At the end these constructions are applied to basic examples of  
noncommutative bundles over a coquasitriangular Hopf algebra. \vspace{0.3cm}

{\bf keywords}:  coquasitriangular Hopf algebra, smash product, noncommutative bundles, noncommutative differential calculus, noncommutative connection
\vspace{0.3cm}

\end{abstract}
\end{titlepage}

%% file: intro.tex
\section*{Introduction}
\addcontentsline{toc}{section}{Introduction}          
The notion of a principal fiber bundle is crucial for geometry and gauge theory.
Because of its importance and its far reaching applications people have 
generalized it to the noncommutative setting in one way or another
(see {\sc Schneider \cite{Sch:PHSAHA}, Pflaum \cite{Pfl:QF,Pfl:QGFB}, Schauenburg \cite{Sch:NDHHGEdRK,Sch:HGEDCNDGPFB},
Brzezinski and Maijd \cite{BrzMaj}}). 
In particular, \cite{Sch:PHSAHA} is based on the fact that the noncommutative
version of an affine principal fiber bundle is a Hopf-Galois extension. Smash
products of a Hopf algebra with a module algebra are important examples of
such extensions. In \cite{Pfl:QF,Pfl:QGFB} smash products (or more generally
crossed products, which are smash products with an additional twisting) 
play the role of local trivializations of fiber bundles
with quantum structure groups. When we view a smash product $A\# H$ of a 
Hopf algebra $H$ with an $H$-module algebra $A$ as a noncommutative principal
fiber bundle, the Hopf algebra $H$ is regarded as the `noncommutative
function space' of the quantum structure group, the algebra $A$ is the 
space of `functions' on the base quantum space, and $A\# H$ itself the 
algebra of `functions' on the quantum fiber bundle.

We will consider noncommutative differential calculi on smash product
`bundles'. Our background for noncommutative differential calculus is based
on the approaches of {\sc Woronowicz} \cite{Wor:DCCMP} and 
{\sc Manin} \cite{Man:NQGQdRC}. In particular, the modules of differential
forms on a quantum space are to be considered as an additional
part of the structure
of such a space, which is not determined by the function space alone, unlike
in the classical case, where there is a canonical functor from commutative
algebras of functions to spaces of differential forms. Whenever the quantum
spaces are endowed with additional structures (like the one of a quantum
group), these have to be compatible with the differentiable structure.

In this spirit we will construct differential calculi on smash products
which are compatible with their algebraic structures. These calculi are
examples for the theories in \cite{BrzMaj,Sch:NDHHGEdRK,Sch:HGEDCNDGPFB},
and they are local differential calculi for the theory in 
\cite{Pfl:QF,Pfl:QGFB}.

Let us summarize the content of our paper.
In the first section we recall the definition of smash products and explain
how they are noncommutative analogues of principal fiber bundles.
In the following sections the concept of first order differential calculi
is introduced, and we construct such a calculus for smash products in a 
canonical way. Later, we lift this first order calculus to a higher order
one.
In section 4 connections on the noncommutative fiber bundles are defined
and their basic structure is studied. We close our article with
an important example.

%% file: sect1.tex
\section{Smash products as noncommutative fiber bundles}
If not 
specified otherwise, 
all algebras and Hopf algebras
are supposed to be defined over the field $k$ and possess a unit $1$.
We denote the comultiplication and counit of a Hopf algebra $H$ by
$\Delta:H\rightarrow H\ot H$ and $\epsilon:H\rightarrow k$. Recall the 
axioms of a Hopf algebra: $\Delta$ is supposed to be an algebra map 
satisfying the coassociativity condition 
$(\Delta\ot 1_H)\Delta=(1_H \ot \Delta)\Delta$, and $\epsilon$ is a counit,
$(\epsilon\ot 1_H)\Delta=(1_H\ot \epsilon)\Delta=1_H$. Furthermore, $H$
has to have an antipode $S:H\rightarrow H$, which, however, we will not 
need explicitly. We will view noncommutative algebras loosely as the algebras
of ``functions'' on ``noncommutative'' or ``quantum spaces''. More precisely
we have a duality (equivalence of categories reversing arrows) between the
categories of quantum spaces and noncommutative algebras. In this picture,
a Hopf algebra is viewed as the algebra of functions on a quantum group,
the comultiplication $\Delta$ corresponding by duality to the multiplication
in the group.

Let $h\in H$ be an element in a Hopf algebra. The image $\Delta(h)$ of $h$
under comultiplication can in general only be written as some finite sum
$\sum h_i'\ot h_i''\in H\ot H$ with $h_i',h_i''\in H$. To simplify
calculations it is customary to use Sweedler's notation
\cite[Sec. 1.2]{Swe:HA}: We write
$\Delta(h)=:\sum_{(h)}h_{(1)}\ot h_{(2)}$ where the individual symbols
$h_{(1)}$ and $h_{(2)}$ cannot be used separately, but only make sense
in bilinear expressions containing both of them. By coassociativity we 
can write 
$(\Delta\ot 1_H)\Delta(h)
  =(1_H\ot\Delta)\Delta(h)=\sum_{(h)}h_{(1)}\ot h_{(2)}\ot h_{(3)}$
etc. A right $H$-comodule is a vector space $V$ with a structure map
or coaction $\rho:V\rightarrow V\ot H$ satisfying coassociativity
$(\rho\ot 1_H)\rho=(1_V\ot \Delta)\rho$ and $(1_V\ot\epsilon)\rho=1_V$.
A left $H$-comodule $V$ has a structure map $\lambda:V\rightarrow H\ot V$
with analogous properties. We write
$\rho(v)=\sum_{(v)}v_{(0)}\ot v_{(1)}$, 
$\lambda(v)=\sum_{(v)}v_{(-1)}\ot v_{(0)}$ etc.
In complicated formulas we will sometimes even omit the summation symbols
in Sweedler's notation. Our general reference for
Hopf algebra theory is {\sc Sweedler} \cite{Swe:HA}.

If $V$ is a right (left) $H$-comodule, we define 
$\rcofix VH:=\{v\in V|v_{(0)}\ot v_{(1)}= v\ot 1\}$ 
(resp. $\lcofix VH:=\{v\in V|v_{(-1)}\ot v_{(0)}=1\ot v\}$)
to be the subset of right (resp. left) coinvariant elements.

Let $R$ be an algebra and a (right) $H$-comodule. We say that $R$ is a
(right) $H$-comodule algebra if the structure map $\rho:R\rightarrow R\ot H$
is an algebra map, that is $\rho(xy)=\rho(x)\rho(y)$ holds in the algebra
$R\ot H$. In the language of quantum spaces this coaction corresponds
to an action of a quantum group $G$ on a noncommutative space $X$. The
subalgebra $\rcofix RH$ of coinvariant elements should be viewed as the 
algebra of functions on the quotient $X/G$ of $X$ under the action of $G$.

An important example of a comodule algebra is the smash product of an
algebra $A$ with a Hopf algebra $H$ that acts on $A$. Let us review
the construction.

\begin{Definition}
Let $H$ be a Hopf algebra with counit $\epsilon$ and $A$ an algebra. An $H$-action on $A$ is given
by a linear map $ H \otimes A \longrightarrow A $ such that
\begin{enumerate}
\item
 $ h \cdot (ab) = \sum \limits_{(h)} \, 
 \left( h_{(1)} \cdot a \right) \, \left( h_{(2)} \cdot b \right)
 , \quad h \in H, \quad  a,b \in A $
\item
 $ 1 \cdot a = a, \quad a \in A$
\item
 $ h \cdot 1 = \epsilon ( h) \, 1, \quad h \in H$
\item
 $ g \cdot ( h \cdot a ) = (g h) \cdot a \quad g,h \in H ,\quad a \in A$
\end{enumerate}
In this case $A$ is also called an $H$-module algebra.
\end{Definition}

\begin{Example} \mbox{ }
 Let $H \otimes A \longrightarrow A$ be given by 
 $ h \otimes a \longmapsto h a = \epsilon (h) a $.
 Then $H$ is said to act trivially on $A$.
\end{Example}
For an important example of a nontrivial action we refer to
section 5.

\begin{Theorem and Definition}
 Let $H$ act on the algebra $A$. Then the mapping 
\beas 
 m & : & \, (A \otimes H) \otimes (A \otimes H) \:  \longrightarrow \: (A \otimes H),
 \\
 & & (a \otimes h , b \otimes g ) \longmapsto \sum \limits_{(h)} \,
 a \left( h_{(1)}\cdot b \right) \otimes h_{(2)} g 
\eeas
 defines a product on the linear space $A \otimes H$, such that $(A \otimes H, m)$
 becomes an algebra with unit $1 \otimes 1$.
 We simply write $A \# H$ for this algebra and call it the smash product
 of $A$ and $H$. 
 Through the mapping
\beas
 \rho &: & \, A \otimes H \: \longrightarrow \: (A \otimes H) \ot H \\
 & & a \otimes h \: \longmapsto \: \sum \limits_{(h)} \, \left( a \ot h_{(1)} \right) 
 \ot h_{(2)} 
\eeas
 $A \# H$ becomes 
 an $H$-right comodule algebra. We sometimes write $a \# h := a \otimes h$.
\end{Theorem and Definition}
\bpro see {\sc Montgomery} \cite{Mon:HATAR}, Lemma 7.1.2, where a more general
   statement is proven.
\epro
Let us fix a Hopf-algebra $H$. If $B$ is another $H$-module algebra and
$f : A \longrightarrow B$ an algebra map such that 
\bdm
f(ha) = h f(a), \quad a \in A 
\edm
we call $f$ an $H$-module algebra map.
Now consider the the mapping
\bdm
 f \# 1_H \, : \, A \# H \longrightarrow B \# H \, , \:
 a \ot h \longmapsto f(a) \otimes h .
\edm
It is easily seen that $f \# 1_H$ preserves the product on $A \# H$, and that
the class of all smash products $A \# H$ together with maps $f \# 1_H$ form a
category,
which we call the category of $H$-smash products.
\begin{Definition}
 The category dual to the one of $H$ smash products is called the
 category of noncommutatice $H$-bundles or quantum $H$-bundles.
\end{Definition}
In the duality between quantum spaces and noncommtuative algebras, the 
(cartesian) product of spaces corresponds to the  tensor product of 
algebras. Thus, the tensor product algebra $A\ot H$ is the direct generalization
of the function algebra on a trivial (i.~e.~product) fiber bundle with
structure quantum group corresponding to $H$ and base quantum space represented
by $A$. This is the special case of Theorem 1.3. that results if we take
the trivial action of $H$ on $A$. Taking a nontrivial action results
in a certain twisting of this ``trivial bundle'' which is natural to the
noncommutative setting. An even more general construction is also well known
(see {\sc Montgomery} \cite{Mon:HATAR}), namely, 
crossed products involving a (weak) action of $H$ on $A$ and a cocycle.
We will concentrate on the smash product since for this particular case
we will
be able to construct natural differential calculi in section 2.

In the sequel we will frequently use the notion of an $^H_H$Hopf$\, ^H_H$-module.
For the convenience of the reader we will give its definition.

Let $H$ be a Hopf-algebra and $A$ be an left $H$-comodule algebra with
coaction $\lambda$.
An left $A$-module $V$ which is at the same time
an left $H$-comodule with structure map $\lambda$ is called an $^H_A$Hopf-module,
if
\be
 \lambda (av) \: = \: \lambda a \, \cdot \, \lambda v, \quad a \in A , 
 \quad v \in V.
\ee
Analogously $^H$Hopf$_A$ modules are defined, and, for $B$ a right 
$H$-comodule algebra, $_B$Hopf$^H$ modules and Hopf$^H_B$-modules.
In particular, all four definitions make sense for $A=B=H$.
If $V$ is a left and right $H$-module and a left and right $H$-comodule
with structure maps $\lambda$ and $\rho$ such that
\bea 
 (a v) \, b & = & a \, (v b) \quad a \in A, \quad b \in B, \quad v \in V \\ 
 ( 1_H \otimes \rho ) \circ \lambda  & = & (\lambda \otimes 1_H) \circ \rho
\eea
holds,
and all the four Hopf module compatibility conditions are satisfied,
then we call $V$ an $^H_H$Hopf$\, ^H_H$-module.

%% file: sect2.tex
\section{First order differential calculi}\label{sec:fodc}
\begin{Definition}
A first order differential calculus for an algebra $A$ is an $A$-bimodule
$M$ together with a derivation $d : A \longrightarrow M$, such that
$M$ is generated by $d(A)$ as an $A$-bimodule.
A morphism between first order differential calculi
$(A,M,d)$ and $(B,N,\delta)$ is given by
a pair $(f, F)$, such that the following conditions 
are fulfilled:
\begin{enumerate}
\item 
 $f : A \longrightarrow B$ is a morphism of algebras.
\item 
 $F : M \longrightarrow N$ is an $A$-bimodule map, where $f$ defines the
 $A$-bimodule structure on $N$.
\item 
 The diagram
\begin{center}
\setsqparms[1`1`1`1;700`500]
\square[A`B`M`N;f`d`\delta`F]
\end{center}
commutes.
\end{enumerate}
 An algebra map  $f : A \longrightarrow B$ is called 
 $($once$)$ 
 differentiable, if a map
 $F : M \longrightarrow N$ exists, such that the pair $(f, F)$ is a morphism 
 of first order differential calculi.
\\
 If $(A,M,d)$ is a differential calculus, we usually write
 $\Omega^1(A)$ instead of $A$ and call it the space of $1$-forms.
A map $F$ as above is uniquely determined by $f$, and we will frequently
denote it by $\Omega^1(f):\Omega^1(A)\rightarrow\Omega^1(B)$.
\end{Definition}
Note that if $\Omega^1(A)$ is a first order differential calculus, then
$ \Omega^1(A)$ is generated by $d(A)$ as a left (or right) $A$-module.
Indeed, every $\omega\in\Omega^1(A)$ has, by definition, the form
$\omega=\sum x_id(y_i)z_i$ for some $x_i,y_i,z_i\in A$, and consequently
$\omega=\sum x_id(y_iz_i)-\sum x_iy_id(z_i)\in A\;d(A)$.
\begin{Theorem}
{\rm (Tensor product of differential calculi)} Assume to be given two 
differential calculi $d : A \longrightarrow M $ and 
$\delta : B \longrightarrow N $ over $k$-algebras $A$ and $B$.
Then $D : A \ot B \longrightarrow M \ot B \oplus A \ot N $, 
$a \ot b \longmapsto da \ot b + a \ot \delta b$ is a differential calculus on
the tensor product algebra $A \ot B$. 
\end{Theorem}
\bpro
Obviously $  A \ot B \longrightarrow M \ot B \oplus A \ot N $ is an $A \ot B$-bimodule
in a natural way, and $D$ a $k$-linear map. $D$ is a derivation by the following equation.
\bea
\lefteqn{
 D( ( a \ot b) \, ( a' \ot b') ) \: = \: D ( a \, a' \ot b \, b' ) \: = \: 
 d( a \, a' ) \ot b \, b' + a \, a' \ot \delta ( b \, b' ) \: = } \nonumber \\
 & = & da \cdot a' \ot b \, b' + a \, a' \ot \delta b \cdot b' \, + \, 
 a \cdot da' \ot b \, b' + a \, a' \ot b \cdot \delta b' \: = \\ 
 & = & D( a \ot b ) \cdot ( a' \ot b' ) \, + \, ( a \ot b) \cdot D( a' \ot b' ),
 \hspace{1cm} a \in A ,  \quad b \in B.
 \nonumber
\eea
 As the image $D(A \ot B)$ generates $ M \ot B \oplus A \ot N $, the theorem 
 is shown.
\epro
If we are given differential calculi $\Omega^1 (A)$ and $\Omega^1 (B)$
on $A$ resp.~$B$ we consider the above defined differential calculus on 
$A \ot B$ as the natural one and will always denote it by 
$\Omega^1 (A \ot B)$.
This choice corresponds to a well known fact in the classical case: the
decomposition of $\Omega^1(A\ot B)$ into a direct sum is the analog of the
decomposition of the cotangent space of a product manifold as the direct
sum of the cotangent spaces of its factors.

If we are given a differential calculus on a Hopf algebra, the differential and Hopf
structure should fit together. What we mean by that precisely is explained in the 
following definition.
\bDef
Let $H$ be a Hopf algebra. A differential calculus $(H, \Omega^1 (H), d)$ is
called a Hopf differential calculus, if the comultiplication
$\Delta : H \longrightarrow H \ot H$ is differentiable.
\eDef
\bPro
\label{comdiff}
\marg{\rm label comdiff}
A differential calculus $(H,\Omega^1(H),d)$ on a Hopf algebra $H$ is a
Hopf differential calculus 
if and only if $\Omega^1 (H) $ is an 
$^H_H$Hopf$\, ^H_H$-module such that $d$ is bicolinear.

That is, a first order differential calculus on a Hopf algebra is a Hopf
differential calculus if and only if it is bicovariant in the sense of
{\sc Woronowicz} \cite{Wor:DCCMP}.
\ePro
\bpro
We have $ \Omega^1 (H \ot H) = \Omega^1 (H) \ot H \oplus H \ot \Omega^1 (H)$. 

If $\Delta$ is differentiable there exists an $H$-bimodule map
$\Delta^1 : \Omega^1 (H) \longrightarrow \Omega^1 (H \ot H)$ such that the
diagram 
\begin{center}
\setsqparms[1`1`1`1;1200`500]
\square[H`H \ot H`\Omega^1 (H)`\Omega^1 (H) \ot H \oplus H \ot \Omega^1 (H);\Delta`d`d \ot 1 + 1 \ot d`\Delta^1]
\end{center}
commutes. Now define $\Delta^1_l = pr_l \circ \Delta^1$, 
$\Delta^1_r = pr_r \circ \Delta^1$, where $pr_l$ (resp.~$pr_r$) are the 
projections
from $\Omega^1 (H \ot H)$ onto $H \ot \Omega^1 (H) $ (resp.~$\Omega^1 (H) \ot H$).
$pr_l$, $pr_r$ and $\Delta^1$ are $H$-bimodule maps, and therefore 
so are 
$\Delta^1_l$ and $\Delta^1_r$.
We have for $h \in H$
\bea
 (1 \ot \Delta^1_l) \, \Delta^1_l (dh) & = & \sum_{(h)} \, (1 \ot \Delta^1_l) 
\, ( h_{(1)} \ot d h_{(2)} )\nonumber\\
&=& \sum_{(h)} \, h_{(1)} \ot h_{(2)} \ot d h_{(3)} 
 \\ 
 (\Delta \ot 1) \, \Delta^1_l (dh) & = & \sum_{(h)} \, (\Delta \ot 1) \,
( h_{(1)} \ot d h_{(2)} )\nonumber\\
&=& \sum_{(h)} \, h_{(1)} \ot h_{(2)} \ot d h_{(3)}, 
\eea
which gives
\bea
(1 \ot \Delta^1_l ) \, \Delta^1_l & = & (\Delta \ot 1) \, \Delta^1_l \, ,
\eea
as $\Omega^1 (H)$ is generated by elements of the form $dh$, $h \in H$.  
The same argument proves the corresponding equation for $\Delta^1_r$
\bea
(\Delta^1_r \ot 1) \, \Delta^1_r & = & (1 \ot \Delta ) \, \Delta^1_r 
\eea
and the equations
\bea
(\epsilon \ot 1) \, \Delta^1_l & = & 1,  
\\ 
(1 \ot \epsilon ) \, \Delta^1_r & = & 1. 
\eea 
The compatibility condition 
\bea
(\Delta^1_l \ot 1 ) \Delta^1_r & = & (1 \ot \Delta^1_r ) \Delta^1_l
\eea
is also easily verified:
\bea
(\Delta^1_l \ot 1) \Delta^1_r (dh) = \sum_{(h)} \, h_{(1)} \ot d h_{(2)} \ot h_{(3)} =
(1 \ot \Delta^1_r) \Delta^1_l (dh)
\eea
for $h\in H$,
and $\Omega^1 (H)$ is an \Hopf\ since the left and right comodule structure maps
are $H$-linear. Furthermore the above
diagram entails $d$ being $H$-colinear, and one direction of the proposition
is shown. 

Let us suppose now 
$\Omega^1 (H)$ is
an $H$-comodule with structure maps
$\Delta^1_l :\Omega^1 (H) \longrightarrow H \ot \Omega^1 (H)$ and
$\Delta^1_r :\Omega^1 (H) \longrightarrow \Omega^1 (H) \ot H$.
Define 
\bmd
 \Delta^1 \, : & \Omega^1 (H) & \longrightarrow & \Omega^1 (H \ot H) \\
 & \alpha & \longmapsto & \Delta^1_r (\alpha ) + \Delta^1_l (\alpha ).
\emd  
As $ \Delta^1_l $ and $ \Delta^1_r $ are $H$-bimodule maps, the same holds for $\Omega^1$.
The commutativity of the above diagram is easily seen:
\beas
  \Delta^1 d (h) & = & \Delta^1_r d (h) + \Delta^1_l d (h) = 
 (d \ot 1 ) \Delta (h) + (1 \ot d)  \Delta (h) \\
 & = & (d \ot 1 + 1 \ot d) \Delta (h).
\eeas
Therefore $(\Delta, \Delta^1)$ is a morphism of differential calculi.
\epro
Assume we are given a differential calculus $\delta : A \longrightarrow 
\Omega^1(A) $
on the algebra $A$ and a bicovariant differential calculus 
$d : H \longrightarrow  \Omega^1(H) $
on the Hopf algebra $H$.
Suppose further that $H$ acts on $A$. The question now is whether one can 
construct canonically a differerential calculus on the smash product $A \# H$.
Recall that a smash product is a twisted version of an ordinary tensor
product, which in turn corresponds to a product manifold or in this case
trivial bundle. We also have a natural notion of a tensor product differential
calculus defined on a tensor product of algebras. It will thus be natural
to construct our differential calculus on the ``bundle'' $A\# H$ as a twisting
in some sense of the natural tensor product calculus on the ``trivial bundle''.
In particular, our $\Omega^1(A\#H)$ will have the same direct sum decomposition
as $\Omega^1(A\ot H)$, only the module structures will have to be modified.
This is possible, if $H$ acts on $\Omega^1(A)$:
%, that means if
\begin{Definition}
Let $A$ be an $H$-module algebra and $d:A\rightarrow\Omega^1(A)$ a first order
differential calculus. We say that $H$ acts on $\Omega^1(A)$, or that
$\Omega^1(A)$ is an $H$-module differential calculus, if there is an
$H$-module structure on $\Omega^1(A)$ such that $d$ is an $H$-module map and
\be
\label{Hdiffact}
h \cdot (a \beta) = \sum_{(h)} \, \left( h_{(1)}\cdot a \right) \, 
\left( h_{(2)} \cdot\beta \right)
\ee
\marg{\rm label Hdiffact}
for  $h \in H$, $a \in A$ and $  \beta \in \Omega^1 (A)$.
\end{Definition}
\begin{Remark}
  \begin{enumerate}
     \item Let $A$ and $\Omega^1(A)$ be as in the definition. An $H$-action
        on $\Omega^1(A)$ is uniquely determined, if it exists. Indeed, 
        any $\omega\in\Omega^1(A)$ is a linear combination of elements
        $ada'$ with $a,a'\in A$, and we have 
        $$h\cdot(ada')=(h_{(1)}\cdot a)(h_{(2)}\cdot da')
          =(h_{(1)}\cdot a)d(h_{(2)}\cdot a)$$
        if $(\ref{Hdiffact})$ holds and $d$ is $H$-linear.
     \item If $\Omega^1(A)$ is an $H$-module differential calculus, then a 
        right handed analog of $(\ref{Hdiffact})$
        also holds: we have
        $$ h\cdot (\beta a)=\sum_{(h)}(h_{(1)}\cdot\beta)(h_{(2)}\cdot a)$$
        for $h\in H,\;a\in A$ and $\beta\in\Omega^1(A)$.
  \end{enumerate}
\end{Remark}
Assume we are given an $H$-module differential calculus on $A$. 
The differential calculus on $A \# H$ is then defined by the following terms:
\bea\label{spdcanf}
\begin{array}{lclr}
\Omega^1( A \# H ) & := & \Omega^1 (A) \ot H \oplus A \ot \Omega^1 (H) &
\mbox{ } \hspace{0.5cm} \label{difff} \mbox{(differential forms)} 
\end{array}
\eea
\bea
\begin{array}{lcclr}
D \, : & \! A \# H & \longrightarrow & \Omega^1 (A \# H) & 
\mbox{ } \hspace{0.5cm} \label{der} \mbox{(derivative)} \\
& a \ot h & \longmapsto & \delta a \ot h + a \ot d h \\
\end{array}
\eea
\bea
\begin{array}{lcl}
 (A \# H) \otimes \Omega^1 (A \# H) & \longrightarrow & \Omega^1 ( A \# H) 
 \mbox{ } \hspace{1.2cm}  \label{lefta} \mbox{(left $A \# H$-action)} \\
 (a \ot h, \alpha \ot g + b \ot \xi) & \longmapsto & 
 a \, (h_{(1)}\cdot \alpha )\ot h_{(2)} g + a \,( h_{(1)}\cdot b )\ot h_{(2)} 
 \xi
\end{array}
\eea
\bea\label{spdcend}
\begin{array}{lcl}
 \Omega^1 (A \# H) \otimes (A \# H) & \longrightarrow & \Omega^1 ( A \# H) 
 \mbox{ } \hspace{1.4cm}  \label{righta} \mbox{(right $A \# H$-action)} \\    
 (\alpha \ot g + b \ot \xi,a \ot h) & \longmapsto &
 \alpha \: (g_{(1)}\cdot a ) \ot g_{(2)} h 
  + b \: (\xi_{(-1)}\cdot a) 
 \ot \xi_{(0)} h.
\end{array}
\eea
It is a standard calculation to show 
that the last two of these maps really are 
actions. We prove
the Leibniz-rule for $D$ by the following equation.
\begin{eqnarray}
  \lefteqn{D((a\# g)(b\# h))}\nonumber\\
  &=& \sum_{(g)} \delta\left(a(g_{(1)}\cdot b)\right)\ot g_{(2)}h
      +a(g_{(1)}\cdot b)\ot d(g_{(2)}h)\nonumber\\
  &=& \sum_{(g)}\left(\delta(a)(g_{(1)}\cdot b)+a\delta(g_{(1)}\cdot b)\right)
      \ot g_{(2)}h\\
    &&  +a(g_{(1)}\cdot b)\ot\left(d(g_{(2)})h+g_{(2)}dh\right)\\
  &=& (\delta a\ot g + a\ot dg)(b\# h)+(a\# g)(\delta b\ot h+b\ot dh)\nonumber\\
  &=& D(a\# g)(b\#h)+(a\# g)D(b\# h)\nonumber
\end{eqnarray}
for $g,h\in H$ and $a,b\in A$.  
$D(A \# H)$ generates $\Omega^1(A \# H)$ as an $A \# H$-bimodule. For the
proof of this it suffices to show that $a \, \delta b \ot h$ and $a \ot g \, dh$
are in the subbimodule of $\Omega^1(A \# H)$ generated by $D(A \# H)$.  
But this is clear from
\begin{eqnarray}
 ( a \ot 1 ) \, D(b \ot 1 ) \, ( 1 \ot h) \: = \:
 (a \ot 1) \, (\delta b \ot 1 ) ( 1 \ot h ) \: = \: a  \delta b \ot h, 
\\ 
 ( a \ot g) \, D ( 1 \ot h) \: = \: (a \ot g) \, ( 1 \ot d h ) \: = \: 
\epsilon ( g_{(1)} ) \, a \ot g_{(2)} \, dh 
 \: = \: a \ot g \, dh. 
\end{eqnarray}
The above considerations prove
 the main part of the next theorem.
\begin{Theorem}
\label{rightcovar}
Let $\delta : A \longrightarrow \Omega^1(A) $ and
$d : H \longrightarrow  \Omega^1(H) $
be differential calculi on the algebra $A$ resp.~Hopf algebra $H$ and let 
$H$ act on $\Omega^1 (A)$. Further denote by 
$ \iota : A \longrightarrow A \# H $ the map $a \longmapsto a \ot 1$.
Then there exists a $_{A \# H}$Hopf$\, _{A \# H}^H$-module 
$\Omega^1 (A \# H)$, a derivation
$D : A \# H \longrightarrow \Omega^1 (A \# H)$ and a monomorphism
$i : \Omega^1 (A) \longrightarrow  \Omega^1 (A \# H)$ such that the diagram
\begin{center}
\setsqparms[1`1`1`1;1200`500]
\square[A`A \# H`\Omega^1 (A)`\Omega^1 (A \#H);\iota`d`\delta`i]
\end{center}
commutes, and $D : A \# H \longrightarrow \Omega^1 (A \# H)$ is a 
differential calculus.
We call it the standard differential calculus on $A \# H $.
\end{Theorem}
\bpro 
Let $D : A \# H \longrightarrow \Omega^1 (A \# H)$ be the differential calculus
on $A \# H$ defined through the
formulas (\ref{difff}) to (\ref{righta}). 
As 
\be
i \circ d \, (a) \: = \: da \ot 1 \: = \: D ( a \ot 1 ) \: = \: D \circ \iota (a),
\ee
the above diagram commutes. $ \Omega^1 (A \# H)$ is a right $H$-comodule with
coaction
\beas
 \rho & : & \Omega^1 (A \# H) \: \longrightarrow \: \Omega^1 (A \# H) \ot  H \\
 & & \rho ( \alpha \ot g  + a \ot \gamma ) \: \longmapsto \: \sum \limits_{(g)} \,
 \alpha \ot g_{(1)}
 \ot g_{(2)} \, +  \, \sum \limits_{(\gamma)} \, a \ot \gamma_{(0)} \ot
 \gamma_{(1)}, \\
 & & a \in A, \quad g \in H, \quad \alpha \in \Omega^1 (A), \quad \gamma \in \Omega^1 (H).
\eeas
By an easy calculation which is left to the reader it is shown that the
right $A \# H$-action and the right $H$-coaction are compatible.
So the theorem is proven.
\epro
Apart from the $H$-comodule structure on $\Omega^1(A\# H)$, the standard
first order differential calculus has an additional structure. To understand
it we formulate an analog of Proposition \ref{comdiff}
for differential calculi on comodule algebras.
As the most important example for the comodule algebra $R$ in the following theorem think
of a smash-product $A \# H$.
\begin{Theorem}
  Let $R$ be a right $H$-comodule algebra and $\Omega^1(R)$, $\Omega^1(H)$
  first order differential calculi such that the comultiplication of
  $H$ is differentiable. Then $ \rho:R\rightarrow R\ot H$ is differentiable
  if and only if 
  the following conditions are satisfied:
\begin{enumerate}
  \item
  $\Omega^1(R)$ is a $_R$Hopf$\,_R^H$-module such that
  $d:R\rightarrow \Omega^1(R)$ is $H$-colinear.
  \item
  There is a map
  $\pi:\Omega^1(R)\rightarrow R\co{H}\Omega^1(H)$ of
  $_R$Hopf$\,^H_R$-modules making the diagram
  \begin{equation}
  \label{pi.d}
  \begin{array}{c}
     \setsqparms[1`1`1`1;800`500]
     \square[R`R\co HH`\Omega^1(R)`R\co H\Omega^1(H);
             \cong`d`R\co{}d`\pi]
  \end{array}
  \end{equation}
  commute.
\end{enumerate}  
\end{Theorem}
Here $V\co{H}W:=\{\sum v_i\ot w_i\in V\ot W|
  \sum v_{i(0)}\ot v_{i(1)}\ot w_i=\sum v_i\ot w_{i(-1)}\ot w_{i(0)}\}$
denotes the cotensor product of a right $H$-comodule $V$ and a left
$H$-comodule $W$. 
Note that $R\co H\Omega^1(H)$ is indeed a $_R$Hopf$\,_R^H$-module with
\begin{eqnarray*}
  \rho(\sum r_i\otimes \gamma_i)
    &=&\sum r_i\ot \gamma_{i(0)}\ot\gamma_{i(1)}\\
  x\sum r_i\ot\gamma_i&=&\sum x_{(0)}r_i\ot x_{(1)}\gamma_i\\
  (\sum r_i\ot \gamma_i)x&=&\sum r_ix_{(0)}\ot \gamma_ix_{(1)}
\end{eqnarray*}
for $\sum r_i\ot \gamma_i\in R\co H\Omega^1(H)$ and $x\in R$.

\bpro 
  We have $\Omega^1(R\ot H)=\Omega^1(R)\ot H\oplus R\ot\Omega^1(H)$ with
  projections $\pr_\ell$ onto $R\ot\Omega^1(H)$ and $\pr_r$ onto
  $\Omega^1(R)\ot H$.

If $\rho$ is differentiable, there exists $\rho^1:\Omega^1(R)\rightarrow
\Omega^1(R\ot H)$ making
\begin{center}
  \setsqparms[1`1`1`1;1200`500]
  \square[R`R\ot H`\Omega^1(R)`
    \left(\Omega^1(R)\ot H\right)\oplus\left(R\ot\Omega^1(H)\right);
    \rho`d`d\ot 1+1\ot d`\rho^1]
\end{center}
commute. Define 
$\rho^1_r:=\pr_r\delta^1:\Omega^1(R)\rightarrow\Omega^1(R)\ot H$
and $\pi_0:=\pr_\ell\rho^1:\Omega^1(R)\rightarrow R\ot \Omega^1(H)$.
As in Proposition \ref{comdiff} 
we conclude that $\rho^1_r$ makes $\Omega^1(R)$ an
$_R$Hopf$\,_R^H$-module with $d$ a colinear map. Moreover, we have,
for $r\in R$,
\begin{eqnarray*}
  (1\ot \Delta^1_\ell)\pi_0d(r)
     &=&(1\ot\Delta^1_\ell(r_{(0)}\ot dr_{(1)})\\
     &=& r_{(0)}\ot r_{(1)}\ot dr_{(2)}\\
     &=& (\rho^1_r\ot 1)(r_{(0)}\ot dr_{(1)})\\
     &=& (\rho^1_r\ot 1)\pi_0dr
\end{eqnarray*}
Since elements of the form $dr$ for $r\in R$ generate $\Omega^1(R)$, this
gives $(\rho^1_r\ot 1)\pi_0=(1\ot\Delta^1_\ell)\pi_0$, and thus the image
of $\pi_0$ lies in $R\co H\Omega^1(H)$.

As in \ref{comdiff}, $\pi$ and $\rho^1_r$ are $R$-$R$-linear maps since 
$\rho^1$ and $\pr_r$, $\pr_\ell$ are.

Conversely, if $\rho^1_r:\Omega^1(R)\rightarrow\Omega^1(R)\ot H$ making
$\Omega^1(R)$ a Hopf module and $d$ colinear, and an $R$-$R$-linear map 
$\pi:\Omega^1(R)\rightarrow R\co H\Omega^1(H)$ making (\ref{pi.d}) commute
exist, we can define $\rho^1$ by
$\rho^1(\omega)=\rho^1_r(\omega)+\pi(\omega)$.
\epro 

\begin{Theorem}
  Let $A$ be an $H$-module algebra, $\Omega^1(A)$ a first order differential
calculus on which $H$ acts, and $\Omega^1(H)$ a bicovariant first order
differential calculus.

Then the comodule algebra structure of $A\# H$ is once differentiable with 
respect to the standard differential calculus $\Omega^1(A\# H)$.
\end{Theorem}
\bpro 
  We already know from \ref{rightcovar} that a suitable right comodule structure
on $\Omega^1(A\# H)$ exists. The map $\pi$, if it exists, is the unique
map sending $xdy$ in $\Omega^1(A\# H)$ to 
$\sum_{(x),(y)}x_{(0)}y_{(0)}\ot x_{(1)}dy_{(1)}$.
We claim that the composition of
$\pr_1:\Omega^1(A\# H)\rightarrow A\ot\Omega^1(H)$ with the canonical
isomorphism
\begin{eqnarray*}
  \beta:A\ot\Omega^1(H)&\rightarrow&(A\#H)\co H\Omega^1(H)\\
    a\ot \omega&\mapsto&a_{(0)}\# \omega_{(-1)}\ot \omega_{(0)}
\end{eqnarray*}
has this property: indeed, for $a,b\in A$ and $g,h\in H$ we have
\begin{eqnarray*}
  \beta\,\pr_1\left((a\# g)d(b\# h)\right)
  &=& \beta\left(a(g_{(1)}\cdot b)\ot g_{(2)}dh\right)\\
  &=& a(g_{(1)}\cdot b)\# g_{(2)}h_{(1)}\ot g_{(3)}d(h_{(2)})\\
  &=& (a\# g)_{(0)}(b\# h)_{(0)}\ot (a\# g)_{(1)}d((b\# h)_{(1)})
\end{eqnarray*}
\epro 

%% file: sect3.tex
\section{Higher order differential calculi}

\begin{Definition}
  A differential graded algebra consists of an algebra 
  $A\plex=\bigoplus_{i\geq 0}A^i$ with product $\wedge$ and a differential
  $d:A\plex\rightarrow A\plex$ such that 
\begin{enumerate}
  \item
  $A^i \wedge A^j \subset A^{i+j}$ for $i,j\geq 0$.
  \item
  $d^2=0$ and $d(A^i)\subset A^{i+1}$ for $i \geq 0$.
  \item
  $d(\omega\wedge\nu)=d\omega\wedge\nu+(-1)^i\omega\wedge d\nu$ for
  $\omega\in A^i$ and $\nu\in A^j$.
\end{enumerate}

  Let $A$ be an algebra. A higher order differential calculus on $A$ is a 
  differential graded algebra $\Omega\plex(A)$ with $\Omega^0(A)=A$
  such that for all $n \in \N$ the vector space $\Omega^n(A)$ is spanned by elements
  of the form $a_0da_1\wedge\ldots\wedge da_n$ for $a_i\in A$.

  Assume that $\Omega\plex(A)$ and $\Omega\plex(B)$ are higher order 
  differential calculi on algebras $A,B$. We call an algebra homomorphism
  $f:A \rightarrow B$ differentiable, if there is a morphism 
  $\Omega\plex(f):\Omega\plex(A)\rightarrow\Omega\plex(B)$ of differential
  graded algebras with $\Omega^0(f)=f$.
  In other words
  $\Omega\plex(f):\bigoplus\Omega^n(A)\rightarrow\bigoplus\Omega^n(B)$ is
  an algebra map defined by a collection of maps 
  $\Omega^n(f):\Omega^n(A)\rightarrow\Omega^n(B)$ such that 
  $d \, \Omega^n(f)=\Omega^{n+1}(f)\,d$. 
\end{Definition}

\bLem
  If $f:A\rightarrow B$ is differentiable with respect to the higher order
  differential calculi $\Omega\plex(A)$ and $\Omega\plex(B)$, then the 
  corresponding $\Omega\plex(f)$ is unique.
\eLem
\bpro 
  Obvious, since $\Omega\plex(A)$ is generated as an algebra by $\Omega^0 (A) = A$
  and $dA$.
\epro

As an example for the concept of differentiability, we again consider
differential calculi on Hopf algebras:
\bLem
  Let $H$ be a bialgebra and $\Omega\plex(H)$ a higher order differential
  calculus on $H$. The comultiplication $\Delta$ of $H$ is differentiable
  if and only if $\Omega\plex(H)$ has the structure of a differential graded
  bialgebra with $\Omega\plex(H)=H$ as bialgebras.
\eLem
For a proof see \cite{Sch:DGHAQGDC}. 

We need some notation for the graded comultiplication
of the differential graded bialgebra $\Omega\plex(H)$. For 
$\omega\in\Omega\plex(H)$ we write
\begin{equation}\label{gradedsweedler1}
 \sum_{[\omega]}\omega_{[1]}\ot\omega_{[2]}:=\Omega\plex(\Delta)(\omega)
\end{equation}
Denoting by 
$p^{ij}:\Omega^{i+j}(H\ot H)\rightarrow\Omega^i(H)\ot\Omega^j(H)$
the projection, we define for $\omega\in\Omega^n(H)$ and $i+j=n$
\begin{equation}\label{gradedsweedler2}
  \sum_{<\omega>}\omega_{<1,i>}\ot\omega_{<2,j>}:=p^{ij}\Omega^n(\Delta)(\omega)
\end{equation}
so that
$$\Omega^n(\Delta)(\omega)
  =\sum_{\zwind{i+j=n}{<\omega>}}\omega_{<1,i>}\ot\omega_{<2,j>}.$$

To simplify notation we sometimes omit the summation symbols
$\sum \limits_{<\omega>}$ and $\sum \limits_{[\omega]}$.
\bDef
Let $A$ be an $H$-module algebra and $\Omega\plex(A)$ a differential calculus.
We say that $H$ acts on $\Omega\plex(A)$, if there are module structures
$H\ot\Omega^n(A)\rightarrow\Omega^n(A)$ making $\Omega\plex(A)$ a module
algebra such that $d:\Omega^n(A)\rightarrow\Omega^{n+1}(A)$ are $H$-linear
maps.
\eDef

\bThD\label{hiodcsp}
Let $A$ be an $H$-module algebra, $\Omega\plex(A)$ a higher order differential
calculus on which $H$ acts, and $\Omega\plex(H)$ a differential 
calculus with respect to which $\Delta:H\rightarrow H\ot H$ is differentiable.
Let us define 
$$
  \Omega^n(A\# H):=\bigoplus_{i=0}^n\Omega^i(A)\#\Omega^{n-i}(H)$$
for all $n\geq 0$, and set
\begin{eqnarray*}
  d(\omega\#\gamma)&=&d\omega\#\gamma+(-1)^i\omega\# d\gamma\\
  (\omega\#\gamma)\wedge(\nu\#\gamma')
  &=&(-1)^{jk}\omega\wedge(\gamma_{(-1)}\cdot\nu)\#\gamma_{(0)}\gamma'
\end{eqnarray*}
for $\omega\in\Omega^i(A)$, $\gamma\in\Omega^j(H)$ and
$\nu\in\Omega^k(A)$. Then $\Omega\plex(A\# H)$ is a differential calculus with respect to which $\rho:A\# H\rightarrow A\# H\ot H$ is differentiable.
\eThD
\bpro
To show the associativity of the product we compute
for $\omega\in\Omega^i(A)$, $\omega'\in\Omega^{i'}(A)$, 
$\omega''\in\Omega^{i''}(A)$, $\gamma\in\Omega^j(H)$,
$\gamma'\in\Omega^{j'}(H)$ and $\gamma''\in\Omega^{j''}(H):$
\begin{eqnarray*}
  \lefteqn{\left((\omega\#\gamma)\wedge(\omega'\#\gamma')\right)
    \wedge(\omega''\#\gamma'')}\\
  &=&(-1)^{ji'}\left(\omega\wedge(\gamma_{(-1)}\cdot\omega')
    \#\gamma_{(0)}\gamma'\right)
    \wedge(\omega''\#\gamma'')\\
  &=&(-1)^{ji'+(j+j')i''}\omega\wedge(\gamma_{(-2)}\cdot\omega')
    \wedge(\gamma_{(-1)}\gamma'_{(-1)}
    \cdot\omega'')\#\gamma_{(0)}\gamma'_{(0)}\gamma''\\
  &=&(-1)^{j(i'+i'')+j'i''}\omega\wedge
   \left(\gamma_{(-1)}\cdot(\omega'\wedge(\gamma'_{(-1)}\cdot\omega'')\right)
   \#\gamma_{(0)}\gamma'_{(0)}\gamma''\\
  &=&(-1)^{j'i''}(\omega\#\gamma)\wedge
    \left(\omega'\wedge(\gamma'_{(-1)}\cdot\omega'')
    \#\gamma'_{(0)}\gamma''\right)\\
  &=&(\omega\#\gamma)\wedge\left((\omega'\#\gamma')\wedge(\omega''\#\gamma'')
    \right).
\end{eqnarray*}

Now we have to prove that 
$\rho:A\# H\rightarrow A\# H\ot H$ is differentiable. Let us define 
$\Omega\plex(\rho):\Omega\plex(A\# H)\rightarrow\Omega\plex(A\# H\ot H)$ by
$\Omega\plex(\rho)(\omega\#\gamma)=\omega\#\gamma_{[1]}\ot \gamma_{[2]}$.
To see that $\Omega\plex(\rho)$ is an algebra map, we compute for
$\omega,\omega',\gamma$ and $\gamma'$ as above:
\begin{eqnarray*}
  \lefteqn{\Omega\plex(\rho)\bigl((\omega\#\gamma)
    \wedge(\omega'\#\gamma')\bigr)}\\
    &=&(-1)^{ji'}\Omega\plex(\rho)
       \bigl(\omega\wedge(\gamma_{(-1)}\cdot\omega')
       \#\gamma_{(0)}\wedge\gamma'\bigr)\\
    &=&(-1)^{ji'}
       \omega\wedge(\gamma_{(-1)}\cdot\omega')
       \#(\gamma_{(0)}\wedge\gamma')_{[1]}
       \ot (\gamma_{(0)}\wedge\gamma')_{[2]}\\
    &=&\sum_{\zwind{k+\ell=j}{k'+\ell'=j'}}(-1)^{ji'+\ell k'}
       \omega\wedge(\gamma_{(-1)}\cdot\omega')
       \#\gamma_{(0)<1,k>}\wedge\gamma'_{<1,k'>}\ot\\
    && \qquad\ot \gamma_{(0)<2,\ell>}\wedge\gamma'_{<2,\ell'>}\\
    &=&\sum_{\zwind{k+\ell=j}{k'+\ell'=j'}}(-1)^{ki'+\ell (i'+k')}
       \omega\wedge(\gamma_{<1,k>(-1)}\cdot\omega')
       \#\gamma_{<1,k>(0)}\wedge\gamma'_{<1,k'>}\ot\\
    &&\qquad \ot \gamma_{<2,\ell>}\wedge\gamma'_{<2,\ell>}\\
    &=&\sum_{\zwind{k+\ell=j}{k'+\ell'=j'}}(-1)^{\ell(i'+k')}
      (\omega\#\gamma_{<1,k>})
     \wedge(\omega'\#\gamma'_{<1,k'>})
     \ot\gamma_{<2,\ell>}\wedge\gamma'_{<2,\ell'>}\\
    &=&(\omega\#\gamma_{[1]}\ot\gamma_{[2]})\wedge
      (\omega'\#\gamma'_{[1]}\ot\gamma'_{[2]}). 
\end{eqnarray*}
Since we have
$d\Omega(\Delta)=\Omega(\Delta)d$ for $\Omega(\Delta):\Omega(H)\rightarrow\Omega(H)\ot\Omega(H)$, it follows that
\begin{eqnarray*}
  \Omega^{i+j+1}(\rho)d(\omega\#\gamma)
    &=& \Omega^{i+j+1}(\rho)(d\omega\#\gamma+(-1)^i\omega\# d\gamma)\\
    &=& d\omega\# \gamma_{[1]}\ot \gamma_{[2]}
        +(-1)^i\omega\# (d\gamma)_{[1]}\ot(d\gamma)_{[2]}\\
    &=& d\omega\#\gamma_{[1]}\ot\gamma_{[2]}
        +(-1)^i\gamma\# d(\gamma_{[1]}\ot\gamma_{[2]})\\
    &=& d(\omega\#\gamma_{[1]}\ot\gamma_{[2]})\\
    &=& d\Omega^{i+j}(\rho)(\omega\#\gamma)
\end{eqnarray*}
for $\omega,\gamma$ as above,
and thus $d\Omega(\rho)=\Omega(\rho)d$.
\epro
Note that the module $\Omega^1(A\# H)$ of differential 1-forms in \ref{hiodcsp} 
coincides with the one defined earlier in section \ref{sec:fodc}.

To describe the structure of the modules of differential forms
on $\Omega(A\# H)$ we introduce a general construction of a Hopf module over
$A\# H$ as a semidirect product:
\bDef
Let $A$ be an $H$-module algebra and $N$ an $A$-$A$-bimodule. An $H$-action 
on $N$ is an $H$-module structure on $N$ such that
$$h\cdot(ana')=\sum_{(h)}(h_{(1)}\cdot a)(h_{(2)}\cdot n)(h_{(3)}\cdot a')$$
holds for all $h\in H,n\in N,a,a'\in A$.

If $N$ is an $A$-$A$-bimodule on which $H$ acts, and $\Gamma$ is a 
$_H^H$Hopf$^H_H$-module, then $N\#\Gamma:=N\ot\Gamma$ is a 
$_{A\# H}$Hopf$^H_{A\# H}$-module with structure maps
\begin{eqnarray*}
  (a\# g)(n\# h)&=& a(g_{(1)}\cdot n)\# g_{(2)}h\\
  (n\# h)(a\# g)&=& n(h_{(1)}\cdot a)\# h_{(2)}g\\
   \rho(n\# h)&=& n\# h_{(1)}\ot h_{(2)}
\end{eqnarray*}
for $g,h\in H,a\in A$ and $n\in N$.
\eDef
Note that with this definition we can write
$$\Omega^n(A\# H)=\bigoplus_{i=1}^n \Omega^i(A)\#\Omega^{n-i}(H)$$
as an $_{A\# H}$Hopf$_{A\# H}^H$-module.

%% file: sect4.tex
\section{Connections}
Let us motivate our concept of a noncommutative connection. 
In the classical case a conncection on a fiber bundle corresponds to a splitting
of the tangent bundle of the fiber bundle into horizontal and vertical 
tangent vectors. 
But a connection also gives rise to a splitting for the space of cotangent
vectors. More precisely we have a natural notion of horizontal differential
$1$-forms on any smooth fiber bundle $P$, i.~e.~we have an exact sequence 
\begin{equation}
\label{shortsequ}
  0 \rightarrow \Omega_h^1 (P) \rightarrow \Omega^1 (P) \rightarrow
  \Omega^1 (P) /\Omega_h^1 (P) \rightarrow 0,
\end{equation}
where $\Omega^1_h (P)$ is the pull-back bundle of $\Omega^1 (M)$ to $P$.
A connection on $P$ now gives a splitting of this sequence (into 
horizontal and vertical parts).
In Theorem \ref{shex} we will give the noncommutative analog of this sequence
and later on define certain notions of a noncommutative connection as 
different kinds of splittings  of this sequence.

Suppose $R$ is an $H$-comodule algebra and $\Omega(H)$, $\Omega(R)$ are 
differential calculi with respect to which $\Delta:H\rightarrow H\ot H$ and
$\rho:R\rightarrow R\ot H$ are differentiable. Then we have seen that 
there is a map $\pi^1:\Omega^1(R)\rightarrow R\co H\Omega^1(H)$ of
$_R$Hopf$^H_R$-modules induced by the $\Omega(H)$-comodule structure of
$\Omega(R)$. Suppose that $A:=\rcofix RH$ is equipped with a differential 
calculus $\Omega(A)$ such that the inclusion map of $A$ into $R$ is 
differentiable. Then the map $\iota:\Omega^1(A)\rightarrow\Omega^1(R)$ together
with the multiplication in $\Omega(R)$ induces maps
\begin{eqnarray*}
  m_\ell:R\ou A\Omega^1(A)&\rightarrow&\Omega^1(R)\\
  r\ot \omega&\mapsto&r\iota(\omega)\\
  \\
  m_r:\Omega^1(A)\ou AR&\rightarrow&\Omega^1(R)\\
  \omega\ot r&\mapsto&\iota(\omega)r
\end{eqnarray*}
Speaking in a somewhat loose language the spaces $R\ou A\Omega^1(A)$ 
and $\Omega^1(A)\ou AR$ both correspond naturally to the space of horizontal forms on $R$. That we have a left and right version is a special feature of
the noncommutative setting.

Let us now investigate the maps $m_{\ell}$ and $m_{r}$ in 
the special case of $R=A\# H$ with the differential
calculus constructed in section \ref{sec:fodc}, eqns. (\ref{spdcanf}) to
(\ref{spdcend}).
\begin{Theorem}\label{shex}
  Let $A$ be a left $H$-module algebra, $\Omega(H)$ a differential calculus
  with respect to which $\Delta:H\rightarrow H\ot H$ is differentiable, and
  $\Omega(A)$ a differential calculus on which $H$ acts. Let 
  $\Omega(A\# H)$ be the smash product differential calculus.
  Then we have the following short exact sequences in $\HMod{}H{A\# H}{A\# H}$:
      $$0\rightarrow(A\# H)\ou A\Omega^1(A)\nabb{m_\ell}\Omega^1(A\# H)
    \nabb{\pi^1}(A\# H)\co H\Omega^1(H)\rightarrow 0$$
  $$0\rightarrow\Omega^1(A)\ou A(A\# H)\nabb{m_r}\Omega^1(A\# H)
    \nabb{\pi^1}(A\# H)\co H\Omega^1(H)\rightarrow 0.$$
Note that these two short exact sequences are noncommutative analogues
of (\ref{shortsequ}).
\end{Theorem} 

The proof of the Theorem will be immediate from the following lemma:
\begin{Lemma}
  We have isomorphisms
  \begin{eqnarray*}
    \alpha_r:\Omega^1(A)\ou A(A\# H)&\rightarrow&\Omega^1(A)\# H\\
    \omega\ot a\# h&\mapsto &\omega a\# h\\
    \alpha_\ell:(A\# H)\ou A\Omega^1(A)&\rightarrow&\Omega^1(A)\# H\\
    a\# h\ot \omega&\mapsto& a(h_{(1)}\cdot \omega)\# h_{(2)}\\
    \beta:A\# \Omega^1(H)&\rightarrow&(A\# H)\co H\Omega^1(H)\\
    a\#\gamma&\mapsto&a\#\gamma_{(-1)}\ot \gamma_{(0)}
  \end{eqnarray*}
  making the diagrams
\begin{center}
    \settriparms[1`1`1;800]
    \qtriangle[\Omega^1(A)\ou A(A\# H)`\Omega^1(A)\# H`\Omega^1(A\# H);
    \alpha_r`m_r`\iota_1]
\end{center}
\begin{center}
    \settriparms[1`1`1;800]
    \qtriangle[(A\# H)\ou A\Omega^1(A)`\Omega^1(A)\# H`\Omega^1(A\# H);
    \alpha_\ell`m_\ell`\iota_1]
\end{center}
\begin{center}
    \settriparms[1`1`-1;800] 
    \ptriangle[\Omega^1(A\# H)`(A\# H)\co H\Omega^1(H)`A\#\Omega^1(H);
      \pi^1`\pr_2`\beta]
\end{center}
  commute.
\end{Lemma}
\begin{Proof}
  To show commutativity of the three triangles 
  we compute
  $$m_r(\omega\ot a\# h)=(\omega\#1)(a\#h)=\omega a\# h
    =\iota_1\alpha_r(\omega\ot a\# h)$$
  
  $$m_\ell(a\# h\ot \omega)=(a\# h)(\omega\# 1)=a(h_{(1)}\cdot\omega)\# h_{(2)}
    = \iota_1\alpha_\ell(a\# h\ot \omega)$$
    
  $$\beta\pr_2(\omega\# h+a\#\gamma)
  =\beta(a\#\gamma)=a\#\gamma_{(-1)}\ot\gamma_{(0)}
  =\pi^1(\omega\# h + a\#\gamma).$$
  
  It is easy to check that $\alpha_r,\alpha_\ell$ and $\beta$ are isomorphisms
  with $\alpha_r^{-1}(\omega\# h)= \omega\ot 1\# h$, 
  $\alpha_\ell^{-1}(\omega\# h)=1\# S^{-1}(h_{(1)})\cdot \omega\ot h_{(2)}$ and
  $\beta^{-1}(\sum a_i\# h_i\ot\gamma_i)=\sum a_i\# \epsilon(h_i)\gamma_i$.
 \end{Proof}
 
Note that the lemma actually shows that the two short exact sequences in 
theorem \ref{shex}  are just the sequence
\begin{equation}
\label{noncommshortsequ}
0\rightarrow \Omega^1(A)\# H\rightarrow\Omega^1(A\# H)
  \rightarrow A\#\Omega^1(H) \rightarrow 0
\end{equation}
in disguise.  Nevertheless
the complicated statement in \ref{shex} is more natural, as it uses only
maps whose existence is a consequence of differentiability assumptions and
not of the specific construction of $\Omega^1(A\# H)$.
In fact, the 
exact sequence above can also be studied in more general cases, 
\cite{Sch:HGEGAFFC,Sch:HGEDCNDGPFB}.

It is obvious that the module of vector fields with respect to some chosen
differential calculus on an algebra should be the dual module of the module
of differential forms. In our noncommutative setting we have to make a choice
with respect to which side we want to dualize our bimodules of differential
forms. We make our choice so that it is natural to write a differential 
form to the left of a vector field it is evaluated on, as it is 
common in differential geometry.

\begin{Definition}
  Let $A$ be an algebra. We denote the functor
  $$\HMod{}{}AA\ni N\mapsto \Hom_{A-}(N,A)\in\HMod{}{}AA$$
  assigning to an $A$-$A$-bimodule the bimodule of left linear maps to $A$
  by $^*(\mbox{---})$.
  
  Let $A$ be an algebra and $\Omega^1(A)$ a first order differential calculus
  on $A$. The module of right vector fields on $A$ $($with respect to 
  the chosen differential calculus$)$ is 
  ${\cal X}_r(A):=\Hom_{A-}(\Omega^1(A),A)=:{^*}\Omega^1(A)$
  
  We shall write $\langle \omega,X\rangle \in A$ for the evaluation of 
  $\omega\in\Omega^1(A)$ on $X\in{\cal X}_r(A)$.
\end{Definition}

We recall the definition of the Lie algebra of a bicovariant first order
differential calculus on a Hopf algebra from {\sc Woronowicz} \cite{Wor:DCCMP}.
To motivate it note that classically the Lie algebra of a Lie group is the
space of left invariant vector fields, which is dual to the space of left
invariant differential one forms.
We use the notation described in \cite{Sch:DGHAQGDC}
\begin{Definition}
  Let $H$ be a Hopf algebra and $\Omega^1(H)$ a first order differential
  calculus on $H$ with respect to which the comultiplication is differentiable.
  Then we denote by  $\h=\Hom_k(\lcofix{\Omega^1(H)}H,k)$ the Lie algebra of $H$
  with respect to $\Omega^1(H)$. Assume that $\lcofix{\Omega^1(H)}H$ is finite 
  dimensional and $H$ has a bijective antipode. Then we equip $\h$ with the
  right comodule structure defined by
  $\langle\gamma,X_{(0)}\rangle X_{(1)}
  =\langle\gamma_{(0)},X\rangle S^{-1}(\gamma_{(1)})$ for all $X\in \h$
  and $\gamma\in\lcofix{\Omega^1(H)}H$. This is the comodule structure of
  $\h$ as the right dual of $\lcofix{\Omega^1(H)}H$ in the category of
  right $H$-comodules, that is, the unique comodule structure for which 
  the evaluation map $\lcofix{\Omega^1(H)}H\ot\h\rightarrow k$ is 
  colinear.
\end{Definition}
In the following we assume $\Omega^1(H)$ being finite, that means 
$\lcofix{\Omega^1(H)}H$ is finite dimensional. Note that the map
$k\rightarrow \h\ot \lcofix{\Omega^1(H)}H$ defined by $1_k \mapsto \sum x_i\ot x^i$, where $x_i$ is a basis of $\h$ and $x^i$ the dual
basis, is also colinear.

\begin{Lemma}\label{leftfree}
  $A\#\Omega^1(H)$ is a free left $A\# H$-module on the vector subspace
  of all $1\#\gamma$ for $\gamma\in\lcofix{\Omega^1(H)}H$; that is, the
  map 
  \begin{eqnarray*}
    A\# H\ot\lcofix{\Omega^1(H)}H&\rightarrow &A\#\Omega^1(H)\\
    a\# h\ot\gamma&\mapsto&(a\# h)(1\# \gamma)=(a\# h\gamma)
  \end{eqnarray*}
  is an isomorphism.
\end{Lemma}
\begin{Proof} 
  By the structure theorem on Hopf modules, {\sc Sweedler} \cite{Swe:HA}, we have
  \begin{eqnarray*}
    H\ot\lcofix{\Omega^1(H)}H&\cong&\Omega^1(H)\\
    h\ot\gamma&\mapsto& h\gamma\\
    \gamma_{(-2)}\ot S(\gamma_{(-1)})\gamma_{(0)}&\mapsfrom&\gamma
  \end{eqnarray*}
  whence an inverse for the map of the lemma is given by
  $$ A\#\Omega^1(H)\ni a\#\gamma
     \mapsto (a\#\gamma_{(-2)})\ot S(\gamma_{(-1)})\gamma_{(0)}
     \in A\# H\ot\lcofix{\Omega^1(H)}H.$$
\end{Proof}

\begin{Definition and Lemma}
  Let $X\in\h$. The fundamental vector field $\overline X$ on $A\# H$ 
  associated to $X$ is the unique right vector field 
  $\overline X\in{\cal X}_r(A\# H)$ vanishing on $\Omega^1(A)\# H$ and 
  satisfying $\langle 1\# \gamma,\overline X\rangle=\langle\gamma,X\rangle 
  1_{A\#H}$
  for all $\gamma\in\lcofix{\Omega^1(H)}H$. We have
  $\langle \omega\# h+a\# \gamma,\overline X\rangle=(a\#\gamma_{(-2)})
  \langle S(\gamma_{(-1)})\gamma_{(0)},X\rangle$ for all 
  $a\in A$, $\omega\in\Omega^1(A)$, $\gamma\in\Omega^1(H)$ and $h\in H$.
\end{Definition and Lemma}
 The first condition in the definition can be interpreted as saying that
 $\overline{X}$ is a vertical vector field. The space of vertical vector
fields is dual (as an $A$-module) to $A\#\Omega^1(H)$, and the second 
condition in the definition describes a natural embedding of $\h$ into this
dual space.

\begin{Proof}
  Since $\Omega^1(A\# H)\cong\Omega^1(A)\# H\oplus A\#\Omega^1(H)$ as
  left $A\# H$-modules, giving a right vector field on $A\# H$ vanishing
  on $\Omega^1(A)\# H$ is equivalent to giving a left $A\# H$-linear map
  $f:A\#\Omega^1(H)\rightarrow A\# H$. By the preceeding Lemma we know
  that $A\#\Omega^1(H)$ is a left $A\# H$-free module over the vector subspace
  of all $1\#\gamma$ for $\gamma\in\lcofix{\Omega^1(H)}H$, whence there
  is a unique such map $f$ with $f(1\#\gamma)=\langle\gamma,X\rangle \cdot 
  1_{A\#H}$
  for $\gamma\in\lcofix{\Omega^1(H)}H$. For the general formula, we compute
  \begin{eqnarray*}
    \langle \omega\# h+a\#\gamma,\overline X\rangle
      &=& \langle a\#\gamma,\overline X\rangle\\ 
      &=& \langle(a\#\gamma_{(-2)})
         (1\# S(\gamma_{(-1)})\gamma_{(0)}),\overline X\rangle\\
     &=& (a\#\gamma_{(-2)})
        \langle(1\# S(\gamma_{(-1)})\gamma_{(0)}),\overline X\rangle\\
     &=& (a\#\gamma_{(-2)})\langle S(\gamma_{(-1)})\gamma_{(0)},X\rangle
   \end{eqnarray*}
   for $a\in A$ and $\gamma\in\Omega^1(H)$.
\end{Proof}  
\begin{Definition}
  Let $A$ be an algebra and $\Omega\plex(A)$ a differential calculus on $A$.
  Let $V$ be a vector space. The space of $V$-valued differential $n$-forms
  on $A$ is $V\ot\Omega^n(A)$. 
\end{Definition}

A $V$-valued differential $n$-form $\omega$ on $A$ is a possibly indecomposable
tensor in $V\ot\Omega^n(A)$. To simplify calculations, we shall nevertheless
frequently write formally $\omega=\omega_V\ot\omega_A$. The usage of this 
symbolic notation is similar to that of Sweedlers notation.

%In particular we can write, in the situation of the definition,
%\begin{eqnarray*}
 % f(\omega)&=&f(\omega_V)\ot\omega_A\\
  %d(\omega)&=&\omega_V\ot d(\omega_A)
%\end{eqnarray*}

\begin{Definition}
  A $V$-valued differential $0$-form on $A$ will be called a $V$-valued function.
  For $v\in V$ we let $\const_v:=v\ot 1_A\in V\ot A$.

  Let $\omega$ be a $V$-valued differential $1$-form and  
  $X\in{\cal X}_r(A)$ a right vector field. Then we put
  $$\skp{\omega}{X}:=\omega_V\skp{\omega_A}{X},$$
  so that $\skp{\omega}{X}$ is a $V$-valued function.
   
  Let $\omega$ be a $V$-valued differential $n$-form.
  If  
  $f:\Omega^n(A)\rightarrow M$ is a map, we abbreviate
  $$f(\omega):=\omega_V\ot f(\omega_A)=(\id_V\ot f)(\omega).$$
  
  If moreover
  %$\omega$ is a $V$-valued differential $n$-form and 
  $\phi\in V^*$, then we put
  $$\skp{\phi}{\omega}:=\skp{\phi}{\omega_V}\omega_A\in\Omega^n(A).$$
\end{Definition}

%\begin{Definition} Let $A$ be an algebra and $\Omega\plex(A)$ a differential
%calculus on $A$. Let $V$ and $W$ be vector spaces. Let 
%$\omega$ and $\nu$ be a $V$-valued $n$-form and a $W$-valued $m$-form
%on $A$,
%respectively. Then
%$$\omega\wedge\nu:=(V\ot %W\ot\nabla^{mn})(V\ot\tau\ot\Omega^n(A))(\omega\ot\nu)$$
%is a $V\ot W$-valued differential $(m+n)$-form on $A$.
%If $X$ is another vector space and $f:V\ot W\rightarrow X$ is a linear map,
%then
%$$f(\omega\ot\nu):=(f\ot\nabla^{nm})(V\ot\tau\ot\Omega^m(A))(\omega\ot\nu)
 % =f(\omega\wedge\nu)$$
%is an $X$-valued differential $(m+n)$-form on $A$.
%\end{Definition}
%With the symbolic notation introduced above we can write
%\begin{eqnarray*}
 % \omega\wedge\nu&=&\omega_V\ot\nu_W\ot\omega_A\wedge\nu_A\\
  %f(\omega\ot\nu)&=&f(\omega_V\ot\nu_W)\ot\omega_A\wedge\nu_A
%\end{eqnarray*}

\begin{Definition}
  Let $A$ be an $H$-comodule algebra and $\Omega\plex(A)$ an $H$-colinear
  differential calculus on $A$. Let $V$ be an $H$-comodule. A $V$-valued
  differential $n$-form $\omega$ is said to be invariant, if we have
  $\omega\in\rcofix{(V\ot\Omega^n(A))}H$, where the $H$-comodule structure 
  on the tensor product is the codiagonal one.
\end{Definition}
Thus $\omega$ is invariant if and only if
$$\omega_{V(0)}\ot\omega_{A(0)}\ot\omega_{V(1)}\omega_{A(1)}=\omega\ot 1$$
Reasoning as in {\sc Schneider} \cite[Lemma 3.1]{Sch:PHSAHA}, this holds iff
$$ \omega_V\ot\omega_{A(0)}\ot\omega_{A\si 1}
  =\omega_{V\si 0}\ot\omega_A\ot S(\omega_{V\si 1}).$$

%\begin{Lemma} Let $A$ be an $H$-comodule algebra, 
%$\Omega\plex(A)$ an 
%$H$-colinear differential calcululs on $A$, let $V$, $W$ and $X$ 
%be $H$-comodules and $f:V\ot W\rightarrow X$ an $H^\op$-colinear map.
%Then for a $V$-valued invariant form $\omega$ and a $W$-valued
%differential form $\nu $ the form $f(\omega\ot\nu)$ is 
%invariant.
%\end{Lemma}
%Indeed we have
%\begin{eqnarray*}
 % \rho(f(\omega\ot\nu))
  %  &=& f(\omega_V\ot\nu_W)\sw 0\ot\omega_{A\si 0}\wedge\nu_{A\si 0}
   %     \ot f(\omega_V\ot\nu_W)\sw 0\omega_{A\si 1}\nu_{A\si 1}\\
    %&=& f(\omega_{V\si 0}\ot\nu_{W\si 0})\ot\omega_{A\si 0}\wedge\nu_{A\si 0}
     %   \ot \nu_{W\si 1}\omega_{V\si 1}\omega_{A\si 1}\nu_{A\si 1}\\
%    &=& f(\omega_V\ot\nu_{W\si 0})\ot\omega_A\wedge\nu_{A\si 0}
 %       \ot \nu_{W\si 1}\nu_{A\si 1}\\
  %  &=& f(\omega_V\ot\nu_W)\ot\omega_A\wedge\nu_A\ot 1
%\end{eqnarray*}

Now we are ready to give the natural definitions for connections and
connection one forms by imitating the classical situation. As announced in
the introduction to this section, a connection, classically a splitting of
the short exact sequence (\ref{shortsequ}), generalizes naturally to a 
splitting of the sequence (\ref{noncommshortsequ}). As for the connection
$1$-forms, we have provided enough notations in the preceding definitions to
be able to just copy literally the classical formula.
\begin{Definition}
  Let $A$ be a left $H$-module algebra, $\Omega(H)$ a differential calculus
  for which $\Delta:H\rightarrow H\ot H$ is differentiable, and 
  $\Omega(A)$ a differential calculus on which $H$ acts.
  \begin{enumerate}
    \item A left connection is a map $c_\ell:(A\# H)\co H\Omega^1(H)\rightarrow
      \Omega^1(A\# H)$ in $\HMod{}H{A\# H}{}$ with $\pi^1 c_\ell=\id$.
    \item A right connection is a map $c_r:(A\# H)\co H\Omega^1(H)\rightarrow
      \Omega^1(A\# H)$ in $\HMod{}H{}{A\# H}$ with $\pi^1 c_r=\id$.
    \item A two-sided connection is a map which is both a left and a right
      connection.  
    \item A connection $1$-form is an element 
      $\phi\in\rcofix{(\h\ot\Omega^1(A\# H))}H$ 
      satisfying
      $\skp{\phi}{\overline X}=\const_X$
      for all $X\in\h$.
  \end{enumerate}
\end{Definition}
\begin{Lemma}
  There exists a canonical two-sided connection on $A\# H$, given by
  $$A\# H\co H\Omega^1(H)
  \nabb{\beta^{-1}}A\#\Omega^1(H)\nabb{\iota}\Omega^1(A\# H),$$
  that is $c(\sum a_i\# h_i\ot \gamma_i)=\sum a_i\epsilon(h_i)\#\gamma_i$.
\end{Lemma}
\begin{Lemma}\label{pivonphi}
  Let $\left\{x_i\right\}$ 
  be a basis of $\h$ and let $\left\{ x^i \right\}$ denote the dual basis of
  $\lcofix{\Omega^1(H)}H$. Then an invariant $\h$-valued $1$-form $\phi$ 
  on $A\# H$
  is a connection $1$-form if and only if 
  $\pi^1(\phi)=\sum x_i \ot 1_{A\#H}\ot x^i$.
\end{Lemma}
\begin{Proof}
  Write $\phi=\varphi+\omega$ with $\varphi\in A\#\Omega^1(H)$ and
  $\omega\in \Omega^1(A)\# H$. We have 
  $\langle\phi,\overline X\rangle=\langle\varphi,\overline X\rangle$, and thus
  it suffices to show that $\varphi=\sum x_i\ot 1\# x^i$ iff 
  $\langle\varphi,\overline X\rangle=X\ot 1_R$ for all $X\in \h$.
  
  If $\varphi=\sum x_i \ot 1\# x^i$, then clearly
  $$\langle \varphi,\overline X\rangle
  =\sum x_i\ot\langle 1\# x^i,\overline X\rangle
  =\sum x_i \ot\langle x^i,X\rangle\cdot 1_R
  =X\ot 1_R.$$
 
  Conversely, assume that $\langle\varphi,\overline X\rangle=X\ot 1$ for all
  $X\in\h$. Write $\pi(\phi)=\sum x_i\ot t_i$ with 
  $t_i\in A\# \Omega^1(H)$.
  By \ref{leftfree}, we can write
  $t_i=\sum r_{ij(0)}\ot r_{ij(1)}x^j$ with $r_{ij}\in A\# H$. It follows that
  \begin{eqnarray*}
    r_{k\ell}
      &=& \sum_j r_{kj}\langle x^j,x_\ell\rangle\\
      &=& \sum_j \langle r_{kj}(1\# x^j),\overline x_\ell \rangle \\
      &=& \sum_{ij} \langle x^k ,x_i\ot 
             \langle r_{ij}(1\# x^j),\overline x_\ell\rangle\rangle\\
      &=& \langle x^k,\langle\varphi,\overline x_\ell\rangle\rangle\\
      &=& \langle x^k,x_\ell\ot 1_{R}\rangle\\
      &=& \rho_{k\ell}\cdot 1_R
  \end{eqnarray*}
  for all $k,\ell$, and thus $\varphi=\sum x_i\ot 1\# x^i$.
\end{Proof}

\begin{Theorem}
  There is a bijection between connection $1$-forms and left connections.
\end{Theorem}
\begin{Proof}
  We identify left connections with left $R$-linear and right $H$-colinear
  maps $A\# \Omega^1(H)\rightarrow \Omega^1(A\# H)$.
  
  Then we can describe the claimed bijection explicitly as follows:
  Given a left connection $c$, we put 
  $$\phi_c=\sum x_i\ot c(1\# x^i) .$$
  Obviously $\pi^1(\phi_c)= \sum x_i \ot 1\ot x^i$, so that $\phi_c$ is a
  connection $1$-form if we verify that it is invariant:
  $\rho(\phi_c)=\sum x_{i(0)}\ot c(1\# x^i_{(0)})\ot x_{i(1)}x^i_{(1)}
    = \sum x_i \ot c(1\# x^i)$
  because $c$ is colinear as well as 
  the map $k\ni 1\mapsto x_i\ot x^i\in\h\ot\h^*$.
  
  Given a connection $1$-form $\phi=\phi_{\h}\ot \phi_R$, 
  put
  $$c_\phi(a\# \gamma)
    =(a\# \gamma_{(-2)})\langle S(\gamma_{(-1)})\gamma_{(0)},\phi\rangle
    =(a\#\gamma_{(-2)})\langle S(\gamma_{(-1)})\gamma_{(0)},\phi_{\h}\rangle
        \phi_R$$
  Let us check that $c_\phi$ is a connection. It is left $R$-linear:
  \begin{eqnarray*}
    c_\phi((x\# h)(a\#\gamma))
      &=& c_\phi(x(h_{(1)}\cdot a)\# h_{(2)}\gamma)\\
      &=& (x(h_{(1)}\cdot a)\# h_{(2)}\gamma_{(-2)})
        \langle S(h_{(3)}\gamma_{(-1)})h_{(4)}\gamma_{(0)},\phi\rangle \\
      &=& (x(h_{(1)}\cdot a)\# h_{(2)}\gamma_{(-2)})
         \langle S(\gamma_{(-1)})\gamma_{(0)},\phi\rangle \\
      &=& (x\# h)(a\#\gamma_{(-1)})
         \langle S(\gamma_{(-1)})\gamma_{(0)},\phi\rangle \\
      &=& (x\# h)c_\phi(a\#\gamma)
  \end{eqnarray*}
  and $H$-colinear:
  \begin{eqnarray*}
    \rho c(a\# \gamma)
      &=& (a\# \gamma_{(-3)})
         \langle S(\gamma_{(-1)})\gamma_{(0)},\phi_{\h}\rangle \phi_{R(0)}
         \ot \gamma_{(-2)}\phi_{R(1)} \\
      &=& (a\# \gamma_{(-3)})
         \langle S(\gamma_{(-1)})\gamma_{(0)},\phi_{\h(0)}\rangle \phi_{R}
         \ot \gamma_{(-2)}S(\phi_{\h(1)}) \\
      &=& (a\# \gamma_{(-4)})
         \langle S(\gamma_{(-1)})\gamma_{(0)},\phi_{\h}\rangle \phi_{R}
         \ot \gamma_{(-3)}S(\gamma_{(-2)})\gamma_{(1)} \\
      &=& (a\# \gamma_{(-2)})
         \langle S(\gamma_{(-1)})\gamma_{(0)},\phi_{\h}\rangle \phi_{R}
         \ot\gamma_{(1)} \\
      &=& c(a\# \gamma_{(0)})\ot \gamma_{1)}.
  \end{eqnarray*}
  Finally
  \begin{eqnarray*}
    \pi^1 c_\phi(a\#\gamma)
      &=& \pi^1((a\#\gamma_{(-2)})
        \langle S(\gamma_{(-1)})\gamma_{(0)},\phi\rangle)\\
      &=& (a\#\gamma_{(-2)})
         \langle S(\gamma_{(-1)})\gamma_{(0)},\pi^1(\phi)\rangle\\
      &=& (a\#\gamma_{(-2)})\langle S(\gamma_{(-1)})\gamma_{(0)},x_i\rangle
          (1\# x^i)\\
      &=& (a\#\gamma_{(-2)})(1\# S(\gamma_{(-1)})\gamma_{(0)})\\
      &=& a\# \gamma. 
  \end{eqnarray*}
  
  In this way we have assigned to every 
  connection a connection $1$-form, and vice versa.
  It remains to verify that the maps thus constructed are
  inverse to each other. But this follows from
  \begin{eqnarray*}
    \phi_{c_\phi}
      &=& \sum x_i\ot c_\phi(1\# x^i)\\
      &=& \sum x_i\ot \langle x^i,\phi_{\h}\rangle\phi_R\\
      &=& \phi_{\h}\ot \phi_R=\phi
  \end{eqnarray*}
  and
  \begin{eqnarray*}
    c_{\phi_c}(a\#\gamma)
      &=& (a\#\gamma_{(-2)})\langle S(\gamma_{(-1)})\gamma_{(0)},\phi_c\rangle\\
      &=& (a\#\gamma_{(-2)})\langle S(\gamma_{(-1)})\gamma_{(0)},x_i\rangle
        c(1\# x^i)\\
      &=& (a\#\gamma_{(-2)})c(1\# S(\gamma_{(-1)})\gamma_{(0)})\\
      &=& c((a\#\gamma_{(-2)})(1\# S(\gamma_{(-1)})\gamma_{(0)}))\\
      &=& c(a\# \gamma).
  \end{eqnarray*}
\end{Proof}

\begin{Remark}
  Surprisingly, there is also a bijection between right connections and
  connection $1$-forms (and consequently also between left and right
  connections). It is given by
   \begin{eqnarray*}
     c_\phi(a\# \gamma)&=&\langle\gamma_{(0)}S^{-1}(\gamma_{(0)},\phi\rangle
       (a\#\gamma_{(-2)}\\
     \phi_c&=&\sum x_i\ot c(1\# x^i).
   \end{eqnarray*}
\end{Remark}
We omit the proof, which is similar to that of the preceding theorem.

\begin{Theorem}
  The set $\cal A$ of connection $1$-forms on $A\# H$ is an affine space
  with translation space isomorphic to the space $\h\ot \Omega^1(A)$
  of $\h$-valued differential $1$-forms on the base quantum space $A$.
\end{Theorem}
\begin{Proof}
  We embed $\h\ot\Omega^1(A)$ into $\h\ot\Omega^1(A\# H)$ via
  $$j:\h\ot\Omega^1(A)\ni X\ot \omega\longmapsto 
     X_{(0)}\ot\omega\# S(X_{(1)})\in\h\ot\Omega^1(A\# H)$$
  and claim that $\cal A$ is an affine space with translation space $\Bi (j)$.

  We have $\Bi(j)\subset\rcofix{\left(\h\ot\Omega^1(A\# H)\right)}{H}$, since
  \begin{eqnarray*}
    \rho(X_{(0)}\ot\omega\# S(X_{(1)}))
      &=&X_{(0)}\ot\omega\# S(X_{(3)})\ot X_{(1)}S(X_{(2)})\\
      &=& X\ot\omega\# S(X_{(1)}).
  \end{eqnarray*}
  Moreover, $\pi^1j(t)=0$ for all $t\in\h\ot\Omega^1(A)$ and thus
  $\phi+j(t)\in{\cal A}$ whenever $\phi\in\cal A$, by \ref{pivonphi}.

  Conversely, let $\phi,\phi'\in\cal A$ and put $\theta=\phi'-\phi$.
  Since $\pi^1(\phi')=\pi^1(\phi)$ by \ref{pivonphi}, we have
  $\pi^1(\theta)=0$, so that $\theta\in\h\ot\Omega^1(A)\#H$. Let us write
  $\theta=\sum x_{ij}\ot \omega_i\# h_{ij}$ with the $\omega_i\in\Omega^1(A)$
  linearly independent, $x_{ij}\in\h$ and $h_{ij}\in H$. Then 
  $\theta\in\rcofix{\left(\h\ot\Omega^1(A\# H)\right)}{H}$ (which follows
  from the fact that both $\phi$ and $\phi' $ are invariant) implies
  $\sum x_{ij}\ot h_{ij}\in\rcofix{(\h\ot H)}{H}$ for all $i$.
  This means $\sum x_{ij}\ot h_{ij}\ot 1
  =\sum x_{ij(0)}\ot h_{ij(1)}\ot x_{ij(1)}h_{ij(2)}$. Now we conclude
  $\sum x_{ij(0)}\ot h_{ij}\ot S(x_{ij(1)})=\sum x_{ij(0)}\ot h_{ij(1)}
  \ot S(x_{ij(1)})x_{ij(2)}h_{ij(2)}
  =\sum x_{ij}\ot h_{ij(1)}\ot h_{ij(2)}$
  which entails 
  \begin{eqnarray*}
    j(\sum x_{ij}\epsilon(h_{ij})\ot\omega_i)
      &=& \sum x_{ij(0)}\ot \omega_i\ot S(x_{ij(1)})\epsilon(h_{ij})\\
      &=&\sum x_{ij}\ot\omega_i\ot h_{ij(2)}\epsilon(h_{ij(1)})\\
      &=&\sum x_{ij}\ot\omega_i\ot h_{ij}=\theta
  \end{eqnarray*}
  and in particular $\theta\in\Bi(j)$.
\end{Proof}

%% file: sect5.tex
\section{A class of examples}

For the following definition of a coquasitriangular or braided bialgebra
see {\sc Larson and Towber} \cite{LarTow:QGQLA} and \cite{Sch:CHAQYBE}.

\begin{Definition}
  A coquasitriangular bialgebra $(H,r)$ consists of a bialgebra $H$ and a 
  convolution invertible map $r:H\ot H\rightarrow k$ satisfying \begin{eqnarray*}
    \sum_{(g),(h)}r(g_{(1)}\ot h_{(1)})g_{(2)}h_{(2)}
      &=&\sum_{(g),(h)}h_{(1)}g_{(1)}r(g_{(2)}\ot h_{(2)})\\
    r(f\ot gh)&=&\sum_{(f)}r(f_{(1)}\ot h)r(f_{(2)}\ot g)\\
    r(fg\ot h)&=&\sum_{(h)}r(f\ot h_{(1)})r(g\ot h_{(2)}).
  \end{eqnarray*}
  for $f,g,h\in H$.
\end{Definition}
Convolution invertible means that there is a map 
$\overline r:H\ot H\rightarrow k$ satisfying
$$\overline r(g_{(1)}\ot h_{(1)})r(g_{(2)}\ot h_{(2)})=\epsilon(g)\epsilon(h)
  =r(g_{(1)}\ot h_{(1)})\overline r(g_{(2)}\ot h_{(2)}.$$

It is an important property of coquasitriangular bialgebras that all their
comodules are also modules in a natural way, cf. \cite{Sch:CHAQYBE}.
\begin{Lemma}\label{monfun}
  Let $(H,r)$ be a coquasitriangular bialgebra. Then there is a monoidal
  (that is, tensor product preserving) functor
  $$^H{\cal M}\rightarrow {_H}{\cal M}$$
  assigning to a left $H$-comodule $V$ the vector space $V$ with the 
  $H$-module structure defined by
  $$h\cdot v=v_{(0)}r(v_{(-1)}\ot h)$$
\end{Lemma}
\begin{Proof}
  The equation really defines an $H$-module, for
  \begin{eqnarray*}
    g\cdot(h\cdot v)
      &=& g\cdot v_{(0)}r(v_{(-1)}\ot h)\\
      &=& v_{(0)}r(v_{(-1)}\ot g)r(v_{(-2)}\ot h)\\
      &=& v_{(0)}r(v_{(-1)}\ot gh)\\
      &=& (gh)\cdot v
  \end{eqnarray*}
  Moreover, the functor thus defined maps the tensor product of two
  comodules $V$ and $W$ to the tensor product of modules:
  \begin{eqnarray*}
    h\cdot(v\ot w)
      &=& v_{(0)}\ot w_{(0)}r(v_{(-1)}w_{(-1)}\ot h)\\
      &=& v_{(0)}r(v_{(-1)}\ot h_{(1)})\ot w_{(0)}r(w_{(-1)}\ot h_{(2)})\\
      &=& h_{(1)}\cdot v\ot h_{(2)}\cdot w
  \end{eqnarray*}
  holds for $v\in V$, $w\in W$ and $h\in H$.
\end{Proof}

Assume we are given a left $H$-comodule algebra $A$. Then lemma \ref{monfun}
entails that we can make $A$ a left $H$-module algebra by setting
$h\cdot a=a_{(0)}r(a_{(-1)}\ot h)$. We can then form the smash product algebra
$A\# H$, whose multiplication is given by
$$(a\# g)(b\# h)=a\;r(b_{(-1)}\ot g_{(1)})b_{(0)}\# g_{(2)}h$$

We can read this as a commutation relation. $A$ and $H$ can be identified
with subalgebras of $A\# H$ via $A\ni a\mapsto a\# 1\in A\# H$ and
$H\ni h\mapsto 1\# h\in A\# H$. Elements $a\in A$ and $h\in H$ satisfy the 
commutation relation
$$ha=r(a_{(-1)}\ot h_{(1)})a_{(0)}h_{(2)}$$
in $A\# H$.

Now assume in addition that we are given a left covariant 
first order differential calculus $\Omega^1(A)$ on $A$ and a bicovariant
first order differential calculus $\Omega^1(H)$ on $H$. Again by lemma 
\ref{monfun}, $H$ acts on $\Omega^1(A)$ by 
$h\cdot \omega=\omega_{(0)}r(\omega_{(-1)}\ot h)$, and we can form the
standard differential calculus $\Omega^1(A\# H)$.
The $A\#H$-$A\#H$-bimodule structure on $\Omega^1(A\# H)$ is described by
\beas
  (a\#g)(\omega\#h)&=&a\;r(\omega_{(-1)}\ot g_{(1)})\omega_{(0)}\#g_{(2)}h\\
 (a\#g)(b\#\gamma)&=&a\;r(b_{(-1)}\ot g_{(1)})\#g_{(2)}\gamma\\
(\omega\#h)(a\#g)&=&\omega \;r(a_{(-1)}\ot h_{(1)})a_{(0)}\#h_{(2)}g\\
(b\#\gamma)(a\#g)&=&b\;r(a_{(-1)}\ot\gamma_{(-1)})a_{(0)}\#\gamma_{(0)}g
\eeas
for $a,b\in A$, $g,h\in H$, $\omega\in\Omega^1(A)$ and $\gamma\in\Omega^1(H)$.
From this we can read off commutation relations between differentials
on $A$ and elements of $H$, and differentials on $H$ and elements in $A$:
\begin{eqnarray*}
  (dh)a&=&r(a_{(-1)}\ot h_{(1)})a_{(0)}(dh_{(2)})\\
  h(da)&=&r(a_{(-1)}\ot h_{(1)})(da_{(0)})h_{(2)}
\end{eqnarray*}
The second equation is not useful if we want to simplify expressions in
$\Omega^1(A\#H)$ by moving all the differentials to the right and all
the algebra elements to the left. We can replace it by
$$(da)h=\overline r(a_{(-1)}\ot h_{(1)})h_{(2)}(da_{(0)})$$
Finally we will become more explicit by recalling examples of coquasitriangular
Hopf algebras (see {\sc Faddeev, Reshetikhin and Takhtajan}
\cite{FadResTak:QLGLA}). 
Let
$(R^{ij}_{k\ell})\in M_{N^2}(k)$ an invertible solution  
of the matrix quantum Yang-Baxter equation.
We write $\hat R^{ij}_{k\ell}=R^{ji}_{k\ell}$.
The FRT algebra $A(R)$ is defined by generators
$(T^i_j|i,j=1,\ldots,N^2)$ and relations 
$\hat R^{ij}_{k\ell}T^k_mT^\ell_n=T^i_kT^j_\ell\hat R^{k\ell}_{mn}$.
It is a coquasitriangular bialgebra with $\Delta(T^i_j)=\sum T^i_k\ot T^k_j$
and with $r:A(R)\ot A(R)\rightarrow k$ defined by 
$r(T^i_j\ot T^k_\ell)=\gamma R^{ik}_{j\ell}$, where $\gamma$ is a nonzero
parameter in $k$. The convolution inverse of $r$ is given by
$\overline r(T^i_j\ot T^k_\ell)=(R^{-1})^{ik}_{j\ell}$.

The quantum groups ${\rm SL}_q(N)$, ${\rm SO}_q(N)$, ${\rm SU}_q(N)$, 
${\rm Sp}_q(N)$, ${\rm SO}_q(3,1)$ are all quotients of FRT algebras 
$A(R)$. To ensure that $r:A(A)\ot A(R)\rightarrow k$ induces a well defined
map $r:H\ot H\rightarrow k$, one is forced to make specific choices of the
parameter $\gamma$. Such choices suitable for the examples listed above 
can be found in {\sc Weixler} \cite[Sec.3.2]{Wei:IQ}.

Now assume that we are given an algebra $A$ generated by $\{x^i\}\subset A$
which is a left $H$-comodule algebra via 
$\rho(x^i)= T^i_j\ot x^j$ (where we have adopted the summation convention).
Then $A$ is at the same time a left $H$-module algebra by
$T^i_j\cdot x^k=\gamma R^{ki}_{\ell j}x^\ell$.
The algebra $A\# H$ is generated by $\{x^i\}\cup\{T^i_j\}$ subject to the 
relations satisfied by the $x^i\in A$ and those satisfied by the $T^i_j\in H$,
plus the additional relations
$$T^i_jx^k=\gamma R^{ki}_{m\ell}x^mT^\ell_j$$.

The usual differential calculi considered on algebras $A$ and $H$ as described
above are such that the {\em left\/} $A$-module $\Omega^1(A)$ and the 
{\em left\/} $H$-module $\Omega^1(H)$ are generated by $\{dx^i\}$ and
$\{dT^i_j\}$, respectively (see {\sc Wess and Zumino} \cite{wesZum:CDCQH}, where
this idea appears first). Then the bimodule $\Omega^1(A\# H)$
is generated by the differentials $dx^i$ and $dT^i_j$ subject to the relations
satisfied by the $dx^i$ and the $x^i$ in $\Omega^1(A)$, those satisfied
by the $dT^i_j$ and the $T^i_j$ in $\Omega^1(H)$, plus the additional 
relations
\begin{eqnarray*}
(dT^i_j)x^k&=&\gamma R^{ki}_{m\ell} x^mT^\ell_j\\
(dx^r)T^s_j&=&\gamma^{-1}(R^{-1})^{rs}_{ki}T^i_j(dx^k)
\end{eqnarray*}
In particular, the {\em left\/} $A\# H$-module $\Omega^1(A\#H)$ is
generated by the differentials $dx^i$ and $dT^i_j$.